\documentclass[10pt]{article}

\def\mathbb{\Bbb}
\usepackage{amsfonts,latexsym,amstext}
\usepackage{amsmath}
\usepackage{amssymb}

\usepackage{comment}

\newtheorem{theorem}{Theorem}[section]
\newtheorem{lemma}[theorem]{Lemma}
\newtheorem{proposition}[theorem]{Proposition}
\newtheorem{definition}{Definition}[section]

\newtheorem{hypothesis}[theorem]{Hypothesis}

\newtheorem{remark}[theorem]{Remark}
\newtheorem{corollary}[theorem]{Corollary}

\makeatletter
\@addtoreset{equation}{section}

\makeatother

\def\qed{{\hfill\hbox{\enspace${ \square}$}} \smallskip}
\def\sqr#1#2{{\vcenter{\vbox{\hrule height .#2pt \hbox{\vrule
 width .#2pt height#1pt \kern#1pt \vrule
width .#2pt} \hrule height .#2pt}}}}
\def\square{\mathchoice\sqr54\sqr54\sqr{4.1}3\sqr{3.5}3}

\def\ds{\begin{displaystyle}}
\def\eds{\end{displaystyle}}
\def\dis{\displaystyle }
\def\<{\langle }
\def\>{\rangle }

\def\R{\mathbb R}

\def\E{\mathbb E}
\def\P{\mathbb P}

\def\F{\mathbb F}

\def\M{\mathbb M}

\def\cala{{\cal A}}
\def\calb{{\cal B}}

\def\calf{{\cal F}}
\def\calg{{\cal G}}

\def\calk{{\cal K}}

\def\calp{{\cal P}}

\def\calu{{\cal U}}

\def\call{{\cal L}}
\def\cals{{\cal S}}

\title{Backward stochastic differential equations
associated to jump Markov processes and applications}

\date{}
\author{Fulvia Confortola, Marco Fuhrman
\\Politecnico di Milano,
Dipartimento di Matematica\\
via Bonardi 9, 20133 Milano, Italy\\
e-mail: fulvia.confortola@polimi.it, marco.fuhrman@polimi.it}

\begin{document}

\maketitle

\begin{abstract}
In this paper we study backward
stochastic differential equations (BSDEs) driven by the compensated random measure
associated to a given pure jump Markov process $X$ on a general state space $K$.
We apply these results to prove well-posedness
of  a class of nonlinear parabolic differential equations
on $K$, that generalize  the Kolmogorov equation
of $X$.  Finally we formulate and solve optimal control problems for
Markov jump processes, relating the value function and the optimal control
law to an appropriate BSDE that also allows to construct probabilistically
the unique solution to the
Hamilton-Jacobi-Bellman equation and to identify it with the value function.
\end{abstract}

\section{Introduction}

In this paper we introduce and solve a class of backward
stochastic differential equations (BSDEs for short) driven by a random measure
associated to a given jump Markov process.
We apply the results to study nonlinear variants of the Kolmogorov equation
of the Markov process and to solve  optimal control problems.

Let us briefly describe our framework.
Our starting point is a pure jump Markov process $X$ on a general
state space $K$.
It is constructed in a usual way starting from a positive measure $A\mapsto\nu(t,x,A)$
on $K$, depending on $t\ge0 $ and $x\in K$ and called rate measure, that specifies
the jump rate function  $\lambda(t,x)=\nu(t,x,K)$ and the
jump measure
$\pi (t,x,A)={\nu(t,x,A)}/{\lambda(t,x)}$.
If the process starts at
time $t$ from $x$ then the distribution of its first
jump time $T_1$  is described by the formula
\begin{equation}\label{jumptimeintro}
\P(T_1>s)=
    \exp\left(-\int_t^s\lambda(r,x)\,dr\right) \,ds,
\end{equation}
and the conditional probability that the process is in $A$ immediately
after a jump at time $T_1=s$ is
$$
  \P(X_{T_1}\in A\,|\,T_1=s ) =
  \pi(s,x,A),
$$
see below for precise statements.
We denote by $\F$ the natural filtration of the process $X$.
Denoting by $T_n$ the jump times of $X$, we consider
the marked point process $(T_n,X_{T_n})$ and the associated
random measure $p(dt\,dy)= \sum_{n} \delta_{(T_n,X_{T_n})} $
on $(0,\infty)\times K$, where $\delta$
denotes the Dirac measure.
 In the markovian case
the dual predictable projection $\tilde p$ of $p$ (shortly, the compensator)
has the following explicit expression
$$
\tilde p(dt\,dy)=\nu(t,X_{t-},dy)\,dt.
$$

In the first part of the paper
we introduce
a class of BSDEs   driven by the compensated random
measure $q(dt\,dy):= p(dt\,dy)-\tilde p(dt\,dy)$ and having the following form
\begin{equation}\label{BSDE-dtintro}
    Y_t+\int_t^T\int_KZ_r(y)\,q(dr\,dy) =
g(X_T) +\int_t^T f(r,X_r,Y_r,Z_r(\cdot))\,dr, \qquad t\in [0,T],
\end{equation}
for given generator $f$ and terminal condition $g$. Here $Y$ is real-valued,
while $Z$ is indexed by $y\in K$, i.e. it is a random field on $K$, with
appropriate measurability conditions, and the generator depends on $Z$ in
a general functional way.
Relying upon the representation theorem for the $\F$-martingales by means
of stochastic integrals with respect to $q$ we can prove several results
on \eqref{BSDE-dtintro}, including existence, uniqueness and continuous
dependence on the data.

In spite of the large literature devoted to random measures
(or equivalently to marked point processes) there are relatively few
results on their connections with BSDEs.
General nonlinear BSDEs driven by the Wiener process were first
solved in \cite{PaPe}. Since then, many generalizations have been
considered where the Wiener process was replaced by more
general processes. Backward  equations driven by random measures
have been studied in \cite{TaLi}, \cite{BaBuPa}, \cite{Roy},
 \cite{KhMaPhZh},
 in view of various applications including stochastic
maximum principle,
partial differential equations of nonlocal type,
quasi-variational inequalities and impulse control. The stochastic
equations addressed in  these papers are driven by
a Wiener process and by a jump process,
but the latter is only considered in the Poisson case.
More general
results on BSDEs driven by random measures
can be found in the paper \cite{Xia}, but  they require
a more involved
formulation; moreover, in contrast to
\cite{TaLi} or
\cite{BaBuPa},
the generator $f$  depends on the process $Z$ in a specific
way (namely as an integral of a Nemytskii operator)
that prevents some of applications that we wish to address, for instance  optimal control
problems.

In this paper $X$ is not defined as a solution of a stochastic equation,
but rather constructed as described above.
While we limit ourselves to the case of a pure jump process $X$,
we can allow  great generality.
Roughly speaking, we can treat all strong Markov jump processes
such that the distribution of holding times admits a rate function $\lambda(t,x)$
as in \eqref{jumptimeintro}:
compare
Remark \ref{processomarkov}-3. The process
$X$ is not required to be time-homogeneous, the holding times
are not necessarily exponentially distributed and can be infinite with
positive probability.
Our main restriction is that the rate measure $\nu$ is uniformly bounded,
which implies that the process $X$ is non explosive.
Our results hold for an arbitrary measurable
state space $K$ (provided one-point sets are measurable) and in particular they
can be directly applied to Markov processes with discrete state
space. We note that
a different formulation of the BSDE is possible for the case
  of finite or countable Markov chains and
   has been studied in \cite{Coh-Ell-1}, \cite{Coh-Ell-2}. In the paper
   \cite{Con-Fuh} we address a class of BSDEs driven by more general random
measures, not necessarily related to a Markov process, but the formulation
is different and more involved, and the corresponding results are less complete.
The results described so far are presented
in section
\ref{sec-backward}, after an introductory  section
 devoted
to notation and preliminaries.

In sections \ref{nlk} and \ref{sec-control} we present two main applications
of the general results on the BSDE \eqref{BSDE-dtintro}.
In section \ref{nlk} we consider a class of parabolic
differential equations on the state space $K$, of the form
\begin{equation}\label{NLKdifferenzialeintro}
\left\{\begin{array}{l}\dis
\partial_t v(t,x)+
\call_t v(t,x)
+ f\Big(t,x,v(t,x), v(t,\cdot)-v(t,x) \Big)=0,\qquad t\in [0,T],\,x\in K,
\\
v(T,x)=g(x),
\end{array}
\right.
\end{equation}
where $\call_t$ denotes the generator of $X$ and $f,g$ are given
functions.  Equation \eqref{NLKdifferenzialeintro} is a non linear variant
of the Kolmogorov equation for the process $X$, the classical equation corresponding
to the case $f=0$.
While it is easy to prove well-posedness of
\eqref{NLKdifferenzialeintro} under boundedness assumptions,
we achieve the purpose of finding a unique solution under much weaker conditions
related to the distribution of the process $X$: see
Theorem \ref{theoremnlk}. We construct the solution $v$
by means of a family of BSDEs
parametrized by
$(t,x)\in [0,T]\times K$:
 \begin{equation}\label{bsdepernlkintro}
   Y_s^{t,x}+\int_s^T\int_KZ_r^{t,x}(y)\,q^{t}(dr\,dy) =
g(X_T) +\int_s^T f(r,X_r,Y_r^{t,x},Z_r^{t,x}(\cdot))\,dr,
\qquad s\in [t,T].
 \end{equation}
 By the  results above there exists a unique solution
$(Y_s^{t,x}, Z_s^{t,x})_{s\in [t,T]}$ and previous estimates on the BSDEs
are used to prove well-posedness of
\eqref{NLKdifferenzialeintro}. As a by-product we also obtain
the  representation formulae
$$
v(t,x)=Y_t^{t,x},
\qquad
    Y_s^{t,x}=v(s,X_s),
\qquad
Z_s^{t,x}(y)= v(s,y)-  v(s,X_{s-}),
$$
which are sometimes called, at least in the diffusive case,
non linear Feynman-Kac formulae.

The second application, that we present in  section
\ref{sec-control} is an optimal control problem.
This is formulated in a  classical way by means of a change
of probability measure, see e.g.  \cite{ElK},  \cite{E}, \cite{B}.
For every fixed $(t,x)\in [0,T]\times K$,
we define a class $\cala^{t}$ of admissible control processes
$u$, and
the cost to be minimized and the corresponding   value function
are
$$
J(t,x,u(\cdot))=\E_u^{t,x}\left[
\int_t^Tl(s,X_s,u_s)\,ds + g(X_T)\right],
\qquad
v(t,x)= \inf_{u(\cdot)\in\cala^{t} }J(t,x,u(\cdot)),
$$
where $g,l$ are given real functions.
Here $\E_u^{t,x}$ denotes the  expectation with respect to another probability $\P_u^{t,x}$,
depending on the control process $u$ and
constructed in such a way that the
compensator  under $\P_u^{t,x}$ equals
$r(s,X_{s-},y,u_s)\, \nu(s,X_{s-},dy)\,ds$ for some function $r$
given in advance as another datum of the control
problem.
The Hamilton-Jacobi-Bellman equation for this problem has
the form \eqref{NLKdifferenzialeintro}
where the generator
is the hamiltonian function
\begin{equation}\label{defhamiltonianintro}
    f( s,x,z(\cdot))=\inf_{u\in U}\left\{
l(s , x,u)+ \int_K z(y) \, (r
(s,x,y,u)-1)\,\nu (s,x,dy)\right\}.
\end{equation}

Optimal control of jump Markov processes is a classical topic in stochastic optimization,
and some the first main results date back several decades: among the earliest
contributions we mention the papers \cite{Bo-Va} and \cite{Bo-Va-Wo}
where, following the dynamic programming approach, the value function of the optimal control problem is characterized as the solution of Hamilton-Jacobi-Bellman, whenever it exists. The results are given under boundedness assumptions on the coefficients.
We refer the reader to the treatise
\cite{Guo-He} for a modern account of the
existing theory; in this book
optimal control problems for continuous time Markov chain  are studied
 in the case of discrete state space and infinite time horizon.

Our approach to this control problem  consists in
introducing  a BSDE of the form \eqref{bsdepernlkintro}, where the generator
is given by \eqref{defhamiltonianintro}.
Under appropriate assumptions and making use of the previous results
we prove that the optimal control
problem has a solution,  that the  value function is the unique solution
to the Hamilton-Jacobi-Bellman  equation
and that the value function and the optimal
control can be represented by means of the solution to the BSDE.
This approach based on BSDEs equations allows to treat
in a unified way a large class of control problems, where the state space is general and the
running and final cost are not necessarily bounded;
moreover it allows to construct probabilistically
a solution of the Hamilton-Jacobi-Bellman equation and to identify it with
the value function.
As in optimal control for diffusive processes (perhaps with the exception
of some recent results) it seems that the approach via BSDEs is limited
to the case when the controlled processes have laws that are all
absolutely continuous with respect to a given, uncontrolled process.
More general cases can be found for instance in  \cite{Guo-He} or,
for more general classes of Markov processes, in \cite{Da-bo}:
see also a more detailed comment
in Remark \ref{confrontocontrollo} below.

We finally mention that
the results of this paper admit several variants and
generalizations: some of them are not included here for reasons of brevity and some
are presently in
preparation. For instance, the
Lipschitz assumptions on the generator of the BSDE can be relaxed,
along the lines of the many results available
in the diffusive case, or extensions to the case of vector-valued
process $Y$ or of random time interval can be considered.

\section{Notations, preliminaries and basic assumptions.}
\label{notations}

\subsection{Jump Markov processes}\label{jumpmarkov}
We recall the definition of a Markov process  as given, for instance,
in \cite{GiSk}.  More precisely we will consider a normal, jump Markov process,
with respect to the natural filtration, with infinite lifetime (i.e. non explosive),
in general non homogeneous in time.

Suppose we are given a measurable space $(K,\calk)$,
a set $\Omega$ and a function
$X:\Omega\times [0,\infty)\to K$.
For every $I\subset [0,\infty)$ we
   denote $\calf_I=\sigma(X_t,\;t\in I)$.
 We suppose that for every
 $t\in [0,\infty)$, $x\in K$ a
 probability $\P^{t,x}$ is given on $(\Omega,\calf_{[t,\infty)})$
 and that the following conditions hold.

\begin{enumerate}
  \item $\calk$ contains all one-point sets.
  $\Delta$ denotes a point not included in $K$.

\item
  $\P^{t,x}(X_t=x)=1$  for every
 $t\in [0,\infty)$, $x\in K$.
 \item  For every $0\le t\le s$ and $A\in \calk$ the function $x\mapsto \P^{t,x}(X_s\in A)$ is $\calk$-measurable.

  \item  For every  $0\le u\le t\le s$, $A\in \calk$  we have
 $\P^{u,x}(X_s\in A|\calf_{[u,t]})=\P^{t,X_t}(X_s\in A)$, $\P^{u,x}$-a.s.

  \item For every $\omega\in\Omega$ and $t\ge 0$ there
exists $\delta>0$ such that $X_s(\omega)=X_t(\omega)$ for $s\in [t,t+\delta]$;
this is equivalent to requiring that all the trajectories of
$X$ have right limits when $K$ is given the
discrete topology (the one where all subsets are open).
  \item
For every $\omega\in\Omega$
the number of jumps of the trajectory $t\mapsto X_t(\omega)$ is finite
on every bounded interval.
\end{enumerate}

 $X$ is called a
 (pure) jump process because of condition 5, and
 a non explosive process because of condition 6.

The class of Markov processes we will consider in this paper
will be described by means of a special form of the joint law $Q$
of the first jump time $T_1$ and the corresponding position $X_{T_1}$.
To proceed formally, we first
fix $t\ge 0$  and $x\in K$ and define   the first jump time
 $T_{ 1} (\omega)=\inf\{s>t\,:\, X_s(\omega)\neq
X_{t}(\omega)\}$,
with the convention that $T_{1}(\omega)=\infty$ if the indicated
set is empty. Clearly, $T_1$ depends on $t$.
Take the extra point $\Delta\notin K$ and define $X_\infty(\omega)=\Delta$
for all $\omega\in\Omega$, so that $X_{T_1}:\Omega\to K\cup\{\Delta\}$
is well defined.
On the extended space
$S:=([0,\infty)\times K)\cup \{(\infty,\Delta)\}$
we consider the smallest $\sigma$-algebra, denoted $\cals$, containing
$\{(\infty,\Delta)\}$ and all sets of $\calb([0,\infty))\otimes\calk$
(here and in the following $\calb(\Lambda)$ denotes
the Borel $\sigma$-algebra of a topological space $\Lambda$).
Then
$(T_1,X_{T_1})$ is a random variable with values in $(S,\cals)$. Its law
under $\P^{t,x}$ will be denoted $Q(t,x,\cdot)$.

We  will assume that $Q$ is constructed
  starting
from a given transition measure from $[0,\infty)\times K$ to $K$,
called  \emph{rate measure} and denoted
$\nu(t,x,A)$, $t\in [0,T]$, $x\in K, A\in \calk$. Thus, we require that
$ A\mapsto \nu(t,x,A)$ is a positive measure on $\calk$
for all $t\in [0,T]$ and $x\in K$, and
$(t,x)\mapsto \nu(t,x,A)$ is $\calb([0,\infty))\otimes\calk$-measurable
for all $A\in \calk$. We also assume
\begin{equation}\label{lambdalimitato}
    \sup_{ t\in[0,T],x\in K}\nu(t,x,K)<\infty, \qquad\qquad
    \nu(t,x,\{x\})=0,\;\;  t\in[0,\infty),x\in K.
\end{equation}
Define
$$
\lambda(t,x)=\nu(t,x,K), \qquad
\pi (t,x,A)=\left\{\begin{array}{ll}
\dis\frac{\nu(t,x,A)}{\lambda(t,x)},& \text{if } \lambda(t,x)>0,
\\
1_A(x),& \text{if } \lambda(t,x)=0.
\end{array}\right.
$$
Therefore $\lambda$ is a nonnegative bounded measurable function
and $\pi$ is a transition probability on $K$
satisfying $\pi(t,x,\{x\})=0$ if $\lambda(t,x)>0$,
and $\pi(t,x,\cdot)=\delta_x$ (the Dirac measure at $x$)
if  $\lambda(t,x)=0$. $\lambda $ is called \emph{jump rate function} and
$\pi$ \emph{jump measure}. Note that we have $\nu(t,x,A)= \lambda(t,x)\pi(t,x,A)$
for all  $t\in[0,T]$, $x\in K$, $A\in \calk$.

Given $\nu$, we will require that for the Markov process $X$ we have,
for  $0\le t\le a< b\le \infty$, $x\in K$, $A\in \calk$,
\begin{equation}\label{jumpkernel}
    Q(t,x, (a,b) \times A)=\int_a^b
    \pi(s,x,A)\;\lambda (s,x)\;
    \exp\left(-\int_t^s\lambda(r,x)\,dr\right) \,ds,
\end{equation}
where
  $Q$ was described above as the law of $(T_1,X_{T_1})$
under $\P^{t,x}$.
Note that \eqref{jumpkernel} completely specifies the probability measure  $Q(t,x,\cdot)$ on
$(S,\cals)$: indeed  simple computations show  that, for $s\ge t$,
\begin{equation}\label{jumpkerneldue}
 \P^{t,x}(T_1\in (s,\infty])  =1- Q(t,x,  (t,s] \times K )
 = \exp\left(-\int_t^s\lambda(r,x)\,dr\right),
\end{equation}
and
we clearly have
\begin{equation}\label{jumpkerneltre}
\begin{array}{lll}
    \P^{t,x}(T_1=\infty)&=&\dis
    Q(t,x, \{(\infty,\Delta)\})=\exp\left(-\int_t^\infty\lambda(r,x)\,dr\right),
    \\
    \P^{t,x}(T_1\le t)  &=&Q(t,x, [0,t] \times K)=0.
    \end{array}
\end{equation}
We may interpret \eqref{jumpkerneldue} as the statement that
$T_1$ has exponential distribution
on $[t,\infty]$ with variable rate $\lambda(r,x)$. Moreover, the probability
$\pi(s,x,\cdot)$ can be interpreted as the conditional probability that $X_{T_1}$
is in $A\in\calk$ given that $T_1=s$; more precisely,
$$
  \P^{t,x}(X_{T_1}\in A, T_1<\infty\,|\,T_1 ) =
  \pi(T_1,x,A)\,1_{T_1<\infty}, \qquad  \P^{t,x}-a.s.
$$

\begin{remark}\label{processomarkov}\begin{em}
\begin{enumerate}
  \item The existence of a jump Markov process satisfying
\eqref{jumpkernel} is a well known fact, see for instance
\cite{GiSk} (Chapter III,  \S 1, Theorems 3 and 4) where it is proved
that $X$ is in addition a strong Markov process.
The nonexplosive character of $X$
 ($T_n\to\infty$) is made possible by our assumption
 (\ref{lambdalimitato}).

 We note that our data only consist initially in a measurable space
$(K,\calk)$ and a transition measure
$\nu$ satisfying \eqref{lambdalimitato}.
The Markov process $(\Omega,X,\P^{t,x})$
can be constructed in an arbitrary way provided
 \eqref{jumpkernel} holds.

  \item In \cite{GiSk} (Chapter III,  \S 1, Theorem 2)  the following is also
  proved: starting from $T_0=t$ define inductively
$
T_{n+1}=\inf\{s>T_n\,:\, X_s\neq
X_{T_{n}}\},
$
with the convention that $T_{n+1} =\infty$ if the indicated
set is empty; then, under the probability
$\P^{t,x}$, the sequence $(T_n,X_{T_n})_{n\ge 0}$ is a discrete-time
Markov process in $(S,\cals)$ with transition kernel $Q$, provided we extend
the definition of $Q$ making the state $(\infty,\Delta)$ absorbing, i.e. we define
$$
 Q(\infty,\Delta, [0,\infty)\times K)=0,\qquad
 Q(\infty,\Delta, \{(\infty,\Delta)\})=1.
 $$
 Note that $(T_n,X_{T_n})_{n\ge 0}$ is time-homogeneous although $X$
 is not, in general.

 This fact allows for a simple description of the process $X$.
 Suppose one starts with  a  discrete-time Markov process $(\tau_n,\xi_n)_{n\ge 0}$
  in $S$ with transition probability kernel $Q$ and a given starting
  point $(t,x)\in [0,\infty)\times K$ (conceptually, trajectories
  of such a process are easy to simulate). One can then define a process $Y$
  in $K$  setting $Y_t=\sum_{n=0}^N\xi_n 1_{[\tau_n, \tau_{n+1}) }(t)$,
  where $N=\sup\{n\ge 0\,:\, \tau_n<\infty\}$. Then  $Y$ has the same law
  as the process $X$ under $\P^{t,x}$.

  \item  We comment on the special form \eqref{jumpkernel} of the kernel $Q$,
   that may seem somehow strange at first sight.
   In \cite{GiSk} (Chapter III,  \S 1) it is proved that for a general jump Markov
   process with the strong Markov property the kernel $Q$ must have the form
$$
Q(t,x, (a,b) \times A)=-\int_a^b
    \pi(s,x,A)\; q(x,t,ds), \qquad 0\le t\le a< b\le \infty, x\in K, A\in \calk,
$$
   where $q(x,t,s)= \P^{t,x}(T_1>s)$ is the survivor function of $T_1$ under
   $\P^{t,x}$ and each $\pi(s,x,\cdot)$ is a suitable probability on $K$. Therefore
   our assumption \eqref{jumpkernel} is basically equivalent to the requirement
   that $q(x,t,\cdot) $ admits a hazard rate function $\lambda(s,x)$ (which turns
out to be  independent of $t$ because of the Markov property). Because of
the clear probabilistic interpretation of $\lambda$ and $\pi$, or equivalently
of $\nu$, we have
preferred to start with the measure $\nu$ as our basic object.

  \item Clearly, the class of processes we consider includes as a very special case
 all the time-homogeneous, nonexplosive, jump Markov processes, which
correspond to the function $\nu$ not depending on $t$.
In this time-homogeneous case the only restriction we retain
is the boundedness  assumption (\ref{lambdalimitato}) on the rate function.

In the time-homogeneous case with $K$ a finite or countable set, the matrix
$(\nu(x,\{y\})_{x,y\in K}$ is the usual matrix of transition rates (or $Q-$matrix)
and $(\pi(x,\{y\})_{x,y\in K}$ is the stochastic transition matrix of the
embedded discrete-time Markov chain.
\end{enumerate}

\end{em}
\end{remark}

\subsection{Marked point processes and the associated martingales}

In this subsection we recall some basic facts following \cite{J}.
In the following we fix
a pair $(t,x)\in [0,\infty) \times K$ and look at the process $X$ under
the probability $\P^{t,x}$.
For every
 $t\ge 0$ we denote  $\F^t$ the filtration $(\calf_{[t,s]})_{s\in [t,\infty)}$.
 We recall that condition 5 above implies that for every
 $t\ge 0$ the filtration $\F^t$ is
 right-continuous (see \cite{B}, Appendix A2,  Theorem T26).

The predictable $\sigma$-algebra
 (respectively, the progressive $\sigma$-algebra)
 on $\Omega\times [t,\infty)$  will be denoted by $\calp^{t}$ (respectively,
 by $Prog^{t}$). The same symbols will also denote the
 restriction to
 $\Omega\times [t,T]$ for  some $T>t$.

For every $t\ge 0$ we define
a sequence $(T_n^t)_{n\ge 0}$ of random variables with values
in $[0,\infty]$ setting
\begin{equation}\label{deftn}
    T_0^t(\omega)=t,\qquad
T_{n+1}^t(\omega)=\inf\{s>T_n^t(\omega)\,:\, X_s(\omega)\neq
X_{T_{n}^t(\omega)}(\omega)\},
\end{equation}
with the convention that $T_{n+1}^t(\omega)=\infty$ if the indicated
set is empty.
Since $X$ is a jump process
we have
$T_{n}^t(\omega)< T^t_{n+1}(\omega)$
 if $T_{n}^t(\omega)<\infty$. Since $X$ is non explosive
we have
$T_{n}^t(\omega)\to\infty$.

For $\omega\in\Omega$ we define a random measure on $((t,\infty)\times K,
\calb((0,\infty))\otimes \calk)$ setting
$$
p^{t}(\omega,C)=\sum_{n\ge 1} 1\Big((T_n^{t}(\omega),
X_{T_n^{t}}(\omega))\in C\Big),
\qquad C\in \calb((t,\infty))\otimes \calk,
$$
where   $1(\ldots)$ is the indicator function.
We also use the notation
$p^{t}(ds\,dy)$ or simply $p(ds\,dy)$. Note that
$$p^{t}((t,s]\times A)=
  \sum_{n\ge 1}1(T^{t}_n\le s)\,1(X_{T_n^{t}}\in A),
  \qquad
 s\ge t,
A\in\calk.
$$

By general results (see \cite{J})
it turns out that
 for every nonnegative $\calp^{t}\otimes \calk$-measurable function
  $H_s(\omega,y)$ defined on $\Omega\times [t,\infty)    \times K$
we have
\begin{equation}\label{defcompensatore}
    \E^{t,x} \int_{t}^\infty \int_K H_s(y)\; p^{t}(ds\,dy)
    =\E^{t,x} \int_{t}^\infty \int _K H_s(y)\;
\nu(s,X_{s},dy)\,ds.
\end{equation}
Note that in this equality we may replace $\nu(s,X_{s},dy)\,ds$ by
$\nu(s,X_{s-},dy)\,ds$.
The random measure $\nu(s,X_{s-},dy)\,ds$ is called the compensator, or the
dual predictable projection, of $p^{t}(ds\,dy)$.

Now fix  $T>t$. If a  real function $H_s(\omega,y)$, defined on
$\Omega\times [t,\infty)    \times K$,
is
 $\calp^{t}\otimes
\calk$-measurable and satisfies
$$
\int_t^T \int_K |H_s(y)| \;\nu(s,X_{s},dy)\,ds <\infty, \qquad
\P^{t,x}-a.s.
$$
then the following stochastic integral can be defined
\begin{equation}\label{defstocint}
  \int_t^s \int_K H_r(y)\; q^{t}(dr\,dy):= \int_t^s \int_K H_r(y)\;
p^{t}(dr\,dy)-\int_t^s \int_K H_r(y) \;\nu(r,X_r,dy)\,dr , \qquad s
\in [t,T],
\end{equation}
   as the difference of ordinary
integrals with respect to $ p^{t}(ds\,dy)$ and $\nu(s,X_{s-}^{t,x},dy)\,ds$. Here and in the
following the symbol $\int_a^b$ is to be understood as an integral
over the interval ${(a,b]}$. We shorten this identity writing
$q^{t}(ds \, dy)= p^{t}(ds\, dy) - \nu(s,X_{s-}^{t,x},dy)\,ds$.
 Note
that
$$
\int_t^s \int_K H_r(y)\; p^{t}(dr\,dy) = \sum_{n\ge 1, T_n^{t}\le s}
H_{T_n^{t}}(X_{T_n^{t}}), \qquad s
\in [t,T],
$$
is always well defined since   $T_n^{t }\to\infty$.

For $m\in [1,\infty)$ we define $\call^{m}(p^{t})$ as
the space of $\calp^{t}\otimes \calk$-measurable real functions
$H_s(\omega,y)$
on $\Omega\times [t,T]\times K$
such that
$$
\E^{t,x}\int_{t}^T \int_K |H_s(y)|^m\; p^{t}(ds\,dy)=
\E^{t,x} \int_{t}^T \int_K  |H_s(y)|^m\;\nu(s,X_{s},dy)\,ds<\infty
$$
(the equality of the integrals follows from \eqref{defcompensatore}).
Given an element $H$ of $\call^{1}(p^{t,x})$, the
stochastic integral \eqref{defstocint}
 turns out to be a  a finite variation martingale.

 We define the space $\call^{1}_{loc}(p^{t})$ as the space
 of those elements $H$ such that
$H\, 1_{(t,S_n]}\in \call^{1}(p^{t})$ for some
increasing sequence of $\F^{t}$-stopping times $S_n$ diverging to $+\infty$.

The key result used in the construction of a solution to BSDEs
 is the integral representation theorem of marked
point process martingales,
which is a counterpart of the well known representation result for
Brownian martingales (see e.g. \cite{Re-Yor} Ch V.3  or  \cite{E}
Thm 12.33).

\begin{theorem}\label{rappresentazione}
Given $(t,x)\in [0,T]\times K$,
let $M$ be an $\F^{t}$-martingale on $[t,T]$ with respect
to $\P^{t,x}$. Then there exists a  process $H \in \call^{1}(p^{t})$ such that
\begin{equation}\label{rapprmart}
    M_s= M_t + \int_t^s \int_K H_r(y)\, q^{t}(dr\; dy),\qquad   s \in [t, T].
\end{equation}

\end{theorem}

\noindent
\textbf{Proof.} When  $M$ is right-continuous the result is well known:
see e.g.  \cite{Da-art},\cite{Da-bo}. The general case reduces to this one
by standard arguments that we only sketch: one first introduces
  the completion $\overline{\F}^{t}$ of the filtration $\F^{t}$
with respect to   $\P^{t,x}$. Then $\overline{\F}^{t}$ satisfies the usual
assumptions, so that $M$ admits a right-continuous modification
   $\bar{M}$, that can be represented as in \eqref{rapprmart}
   by means of
   a process $\bar{H} \in \call^{1}(p^{t})$ and $\overline{\calp}^t\otimes \calk$-measurable,
   where
 $\overline{\calp}^t$ denotes the $\overline{\F}^{t}$-predictable $\sigma$-field.
By monotone class arguments, starting from a simple set of generators
of $\overline{\calp}^t\otimes \calk$, one finally proves that
$\bar{H}$ has a modification $H$ such that \eqref{rapprmart} holds.
\qed

Let us define the generator of the Markov process $X$ setting
$$
\call_t \psi(x)=  \int_K\Big( \psi( y)-\psi( x)  \Big)\,
\nu(t,x,dy), \qquad t\ge 0, x\in K,
$$
for every measurable function $\psi:K\to \R$ for which the integral
is well defined.

We recall the Ito formula for the process $X$, see e.g.
\cite{Da-bo} or
\cite{I-W}.
 Suppose
$0\le t <T$ and let
$v:[t,T]\times K\to \R$ be a
 measurable function such that
 \begin{enumerate}
 \item $s\mapsto v(s,x)$ is absolutely continuous
for every $x\in K$, with time derivative denoted $\partial_sv(s,x)$;
\item
 $\{v(s,y)-v(s,X_{s-}),\; s\in [t,T],y\in K\}$
 belongs to $\call^{1}_{loc}(p^{t})$;
 \end{enumerate}
then, $\P^{t,x}$-a.s.
\begin{equation}\label{Itoformula}
\begin{array}{lll}
v(s,X_s) &=&\dis
v(t,x)+\int_t^s \Big(\partial_r v(r,X_r) +\call_r v(r,X_r)\Big)\, dr
\\&&\dis +
\int_t^s \int_K \Big(v(r,y)-v(r,X_{r-})\Big)\, q^{t}(dr\,dy),
\qquad s\in [t,T]
\end{array}
 \end{equation}
 where the stochastic integral is a local martingale.
In differential notation:
 $$
dv(s,X_s ) = \partial_s v(s,X_s )\,ds+\call_s v(s,X_s^{t,x})\, ds+
 \int_K \Big(v(s,y)-v(s,X_{s-})\Big)\, q^{t}(ds\,dy).
 $$

\section{The backward equation} \label{sec-backward}

Let us assume that $\nu$ is a
a transition measure on $K$ satisfying
\eqref{lambdalimitato}. $X$ denotes the Markov process constructed in
section
\ref{notations}, satisfying conditions 1-6
in subsection \ref{jumpmarkov}  as well as \eqref{jumpkernel}.

Throughout this section
we fix a deterministic terminal time $T > 0$ and a pair
$(t,x)\in [0,T]\times K$. We look at all processes  under
the probability $\P^{t,x}$.
 In the following, especially in the proofs, we will  omit the superscript $t$
and write $\F$, $\calp$, $Prog$, $T_n$, $p(ds\,dy)$, $q(ds\,dy)$,
$\mathcal{L}^2(p)$
instead of $T_n^t$, $\F^t$, $\calp^t$, $Prog^t$, $p^{t}(ds\,dy)$, $q^{t}(ds\,dy)$,
$\mathcal{L}^2(p^{t})$.

We are interested in studying the following family of backward equations
parametrized by $(t,x)$: $\P^{t,x}$-a.s.
\begin{equation}\label{BSDE-dt}
    Y_s+\int_s^T\int_KZ_r(y)\,q(dr\,dy) =
g(X_T) +\int_s^T f(r,X_r,Y_r,Z_r(\cdot))\,dr, \qquad s\in [t,T],
\end{equation}
under the following assumptions on the data $f$ and $g$:

\begin{hypothesis}\label{hyp:BSDE-dt}
\begin{enumerate}
\item The final condition $g:K\to\R$ is $\calk$-measurable and $\E^{t,x} |g(X_T)|^2 < \infty$.

\item The generator $f$ is such that
\begin{itemize}
 \item[(i)] for every  $s \in [0,T]$, $x\in K$,
 $r \in \R$,  $f (s,x,r, \cdot)$ is a mapping $ L^2(K,\calk, \nu(s,x,dy))\to \R$;
 \item[(ii)]  for every bounded and $\calk$-measurable function $ z:K \rightarrow \R$,
the mapping
\begin{equation}\label{misurabilitaf}
  (s,x,r)\mapsto f(s,x,r, z(\cdot))
\end{equation}
is $\calb([0,T])\otimes \calk \otimes \calb(\R)$-measurable;
\item[(iii)] there exist $L \ge 0$, $L' \ge 0$
such that for every $s \in [0,T]$, $x\in K$,
$ r, r' \in \R$, $z,z' \in L^2(K,\mathcal{K}, \nu(s,x,dy) )$,
\begin{equation}\label{flipschitz}
|f(s,x,r,z(\cdot))-f (s,x,r',z'(\cdot))|\leq L'|r-r'|+ L \left(\int_K |z(y)-z'(y)|^2\nu(s,x,dy)\right)^{1/2}.
\end{equation}

\item[(iv)] $\dis\E^{t,x} \int_t^T |f(s,X_s,0,0) |^2 ds <\infty$.
\end{itemize}
\end{enumerate}
\end{hypothesis}

In order to study the backward equation
\eqref{BSDE-dt} we need to check the following
measurability property of $ f(s,X_s,Y_s,Z_s(\cdot))$.

\begin{lemma} Let $f$ be a generator satisfying assumptions
$(i)$, $(ii)$ and $(iii)$.

If $Z \in \mathcal{L}^2(p^{t})$, then
the mapping
\begin{equation}\label{generatoremis}
  (\omega,s,y)\mapsto f(s,X_{s-}(\omega),y, Z_s(\omega,\cdot))
\end{equation}
is $\mathcal{P}^{t} \otimes \calb(\R)$-measurable.

If, in addition,  $Y$ is a
 $Prog^{t}$-measurable process, then $$(\omega,s)\mapsto f(s,X_{s-}(\omega),Y_s(\omega), Z_s(\omega,\cdot))$$ is $Prog^{t}$-measurable.
\end{lemma}

\textbf{Proof.} It is enough to prove the required measurability of the mapping
\eqref{generatoremis}, since the other statement of the lemma follows by composition.

Let $B(K)$ denote the space of  $\calk$-measurable and bounded maps $z: K \rightarrow \R$,
endowed with the supremum norm and the corresponding Borel $\sigma$-algebra
$\calb(B(K))$. Note that $B(K)\subset L^2(K,\calk, \nu(s,x,dy))$ by \eqref{lambdalimitato}.
Consider the restriction of the generator $f$ to $ [0,T] \times K \times \R \times  B(K)$.
By $(ii)$ we have that for all $z \in B(K)$ the function $f(\cdot, \cdot, \cdot,z)$ is $\calb([0,T])\otimes \calk \otimes \calb(\R)$-measurable. Moreover, by \eqref{flipschitz}
 and \eqref{lambdalimitato} it follows that for all $(s,x,y) \in  [0,T] \times K \times \R$ $f(s,x,y,\cdot)$ is continuous. This means that the mapping $f: [0,T] \times K \times \R \times  B(K) \rightarrow \R$ is a Carath\'{e}odory function so, in particular, it is $\calb([0,T])
\otimes \calk \otimes \calb(\R) \otimes  \calb(B(K))$-measurable.

 Now let $Z$ be a bounded $\calp^{t}\otimes \calk$-measurable real function.
 Then, for all $s \in[t,T]$ and $\omega \in \Omega$, $Z_s(\omega, \cdot)$ belongs to $ B(K)$ and,
 using a monotone class argument, it easy to verify that
 the map $(s, \omega) \mapsto Z_s(\omega, \cdot)$ is measurable with respect to $\calp^{t}$ and $\calb(B(K))$.
By composition it follows that
 the mapping $$(\omega,s,y)\mapsto f(s,X_{s-}(\omega),y, Z_s(\omega,\cdot)),$$
  is $\calp^{t}\otimes \calb(\R)$-measurable.

Finally, for general $Z \in \mathcal{L}^2(p^{t})$,
 thanks to the
Lipschitz condition $(iii)$, it is possible to write $$f(t,X_{t-}(\omega),y, Z_t(\omega))
= \lim_{n \rightarrow \infty} f(t,X_{t^-}(\omega), y, Z_t^n(\omega)),
$$
where $Z^n$ is a sequence of bounded and $\calp^{t}\otimes \calk$-measurable real functions converging to $Z$ in $\mathcal{L}^2(p^{t})$. The required measurability follows.
 \qed

We introduce the space $\M^{t,x}$
of the processes $(Y,Z)$ on $[t,T]$
such that $Y$ is real-valued and $Prog^{t}$-measurable,
$Z : \Omega\times [t, T] \times K \rightarrow \R$ is $\mathcal{P}^{t}\otimes \mathcal{K}$-measurable and
$$||(Y,Z)||_{\M^{t,x}}^2
: =\E^{t,x}\int_t^T |Y_s|^2ds
+ \E^{t,x} \int_t^T  \int_K |Z_s(y)|^2 \nu(s,X_{s},dy)\,ds  < \infty.$$

The space $\M^{t,x}$, endowed with this norm, is Banach
space, provided we identify pairs of processes whose difference has norm zero.

\begin{lemma} \label{BSDEcasobanale-dt}
Suppose that $f:\Omega\times [t,T]\to \R$ is  $Prog^{t}$-measurable,
$\xi:\Omega\to \R$ is $\calf_{[t,T]}$-measurable and
$$
\E^{t,x} |{\xi}|^2 +
\E^{t,x} \int_t^T |{f}_s|^2 ds<\infty.
$$
Then there exists a unique pair $(Y,Z)$ in $\M^{t,x}$ solution to
the BSDE
\begin{equation}\label{BSDElineare-dt}
  Y_s+\int_s^T \int_K Z_r(y) \, q^{t}(dr \;dy)=
   \xi +\int_s^Tf_r  \, dr,\qquad s\in [t,T].
\end{equation}
Moreover for all $ \beta \in \R$ we have
\begin{equation}\label{identitaenergia-dt}
\begin{array}{l}\dis
\E^{t,x} e^{\beta s}| {Y}_s|^2 +\beta\,\E^{t,x}\int_s^T  e^{\beta r}|{Y}_r|^2 dr
+ \E^{t,x}\int_s^T \int_K e^{\beta r}|{Z}_r(y)|^2
\nu(r,X_r,dy)\,dr
\\
\dis\qquad
=\E^{t,x} e^{\beta T}|{\xi}|^2
  +2\E^{t,x} \int_s^T  e^{\beta r} {Y}_{r} \,f_r \,dr,
  \qquad s\in [t,T],
\end{array}
\end{equation}
and
\begin{equation}\label{stimaprlineare-dt}
 \E^{t,x}\int_t^T   | {Y}_r|^2 dr
+ \E^{t,x}\int_t^T \int_K | {Z}_r(y)|^2
\nu(r,X_r,dy)\,dr \le 8\,
\E^{t,x}\, |{\xi}|^2 + 8(T+1) \E^{t,x}\,\left[\int_t^T|f_r|^2  \, dr\right].
\end{equation}
\end{lemma}

\noindent{\bf Proof.}
To simplify notation we will drop the superscripts $t,x$
and we write the proof in the case $t=0$.

Uniqueness follows immediately using the linearity of \eqref{BSDElineare-dt} and taking the
conditional expectation given $\calf_{[0,s]}$.

Assuming that $(Y,Z)\in \mathbb{M}$ is a solution, we first prove the identity
\eqref{identitaenergia-dt}.
From the Ito formula applied to $ e^{\beta s}| {Y}_s|^2$ it follows that
$$d( e^{\beta s}| {Y}_s|^2)= \beta  e^{\beta s}| {Y}_s|^2  \, ds + 2 e^{\beta s} {Y}_{s-} dY_s + e^{\beta s}|\Delta {Y}_s|^2.$$

So integrating on $[s,T]$
\begin{eqnarray}
\nonumber e^{\beta s}| {Y}_s|^2 & = & - \int_s^T \beta  e^{\beta r}| {Y}_r|^2  \, dr - 2 \int_s^T \int_K e^{\beta r} {Y}_{r-}  {Z}_r(y) q(dr\,dy) - \sum_{s<r\le T} e^{\beta r}|\Delta  {Y}_r|^2\\
\label{eq:int-Itobis} &  & +e^{\beta T}| {\xi}|^2
+2\int_s^T  e^{\beta r} {Y}_{r}  \,f_r \, dr.
\end{eqnarray}
The process $\int_0^s  \int_K e^{\beta r} {Y}_{r-}  Z_r(y) q(dr\,dy)$
is a martingale, because the integrand process $ e^{\beta r} {Y}_{r-}   {Z}_r(y)$ is in $\call^1(p)$:
 in fact from the Young inequality and \eqref{lambdalimitato} we get
$$
\begin{array}{l}\dis
             \E \int_0^T \int_K  e^{\beta r}| {Y}_{r-}| |  {Z}_r (y)| \nu(r,X_r,dy) \, dr
              \\\dis\qquad
     \le \frac{1}{2}\E \int_0^T \int_K e^{\beta r}| {Y}_{r-}|^2 \nu(r,X_r,dy) \, dr +
      \frac{1}{2}\E  \int_0^T \int_K  e^{\beta r}|  {Z}_r (y)|^2 \nu(r,X_r,dy) \, dr
       \\\dis\qquad
      \le \sup_{t,x}\nu(t,x,K)\,\frac{e^{\beta T}}{2}\E  \int_0^T | {Y}_{r}|^2  \, dr +
      \frac{e^{\beta T}}{2}\E  \int_0^T  | {Z}_r (y)|^2 \nu(r,X_r,dy) \, dr < \infty.
             \end{array}
$$
Moreover we have
\begin{eqnarray*}
\sum_{0<r\le s} e^{\beta r}|\Delta  {Y}_r|^2 & = &
\int_0^t \int_K e^{\beta r}| {Z}_r(y)|^2 \, p(dr\,dy)
  \\
 & = &\int_0^s \int_K e^{\beta r}| {Z}_r(y)|^2 \,q(dr \,dy)+
  \int_0^s \int_K e^{\beta r}| {Z}_r(y)|^2 \nu(r,X_r,dy) \, dr,
 \end{eqnarray*}
where the stochastic integral with respect to $q$ is a martingale.
Taking the expectation in \eqref{eq:int-Itobis} we obtain
\eqref{identitaenergia-dt}.

We now pass to the proof of existence.
The solution $(Y,Z)$ is defined by considering
 the martingale $M_s=\E^{\calf_{[0,s]}}[\xi + \int_0^T f_r \,
dr ]$. By the martingale representation theorem
\ref{rappresentazione}, there exists a  process $Z \in
\call^{1}(p)$ such that
\begin{equation*} M_s= M_0+ \int_0^s \int_K Z_r(y) \;q(dy \,dr) , \qquad s \in [0,  T].
 \end{equation*}
 Define the process $Y$ by the formula
\begin{equation*}Y_s = M_s - \int_0^s f_r \; dr, \qquad s \in [0,  T].
 \end{equation*}
Noting that  $Y_T= \xi$, we easily deduce that the equation \eqref{BSDElineare-dt} is satisfied.

It remains to show that $(Y,Z)\in \M$.
Taking the
conditional expectation, it follows from \eqref{BSDElineare-dt} that
$Y_s= \E^{\calf_{[0,s]}}[\xi + \int_s^T f_r \, dr ]$ so that we obtain
\begin{equation}
|Y_s|^2  \le
2 |\E^{\calf_{[0,s]}}\xi|^2
+2\left|\E^{\calf_{[0,s]}}\int_s^T f_r  \, dr\right|^2\\
 \le
 2\E^{\calf_{[0,s]}} \left[|\xi|^2 +  T\int_0^T|f_r|^2  \, dr\right].
 \label{martingausil}
 \end{equation}
Denoting by $m_s$ the right-hand side of \eqref{martingausil}, we see that $m$ is a martingale
by the assumptions of the lemma. In particular, for every stopping time $S$ with values in $[0,T]$,
we have
\begin{equation}\label{uniftempiarresto}
\E|Y_S|^2\le \E\, m_S=\E\,m_T<\infty
\end{equation}
by the optional stopping theorem. Next we define the increasing sequence of stopping times
$$
S_n=\inf\{s\in [0,T]\,:\,
\int_0^s |Y_r|^2dr +
 \int_0^s  \int_K |Z_r(y)|^2 \nu(r, X_r,dy)\, dr >n\},
 $$
with the  convention $\inf \emptyset =T$. Computing the It\^{o} differential
$d( |{Y}_s|^2)$
on the interval $[0,S_n]$ and proceeding as before we deduce
$$
\E\int_0^{S_n}   | {Y}_r|^2 \,dr
+ \E\int_0^{S_n} \int_K | {Z}_r(y)|^2 \nu(r, X_r,dy)\, dr
\le\E \,|Y_{S_n}|^2
  +2\E \int_0^{S_n}  {Y}_{r}  f_r\, dr.
$$
Using the inequalities
$2{Y}_{r}  f_r\le (1/2)|Y_r|^2 + 2|f_r|^2$
and \eqref{uniftempiarresto} (with $S=S_n$) we find the following estimates
\begin{eqnarray}\label{stimasonesseenne}
   \E\int_0^{S_n}   | {Y}_r|^2 \,dr
 &\le &4\,
\E |\xi|^2 +  4(T+1)\E\int_0^T|f_r|^2  \, dr
\\\label{stimasonesseennebis}
 \E\int_0^{S_n} \int_K | {Z}_r(y)|^2 \nu(r, X_r,dy)\, dr
 &\le&
4\E |\xi|^2 +  4(T+1)\E\int_0^T|f_r|^2  \, dr.
\end{eqnarray}
Setting $S=\lim_nS_n$ we obtain
$$
\int_0^{S } | {Y}_r|^2 \,dr
+  \int_0^{S } \int_K | {Z}_r(y)|^2 \nu(r, X_r,dy)\, dr <\infty,\qquad \P-a.s.
$$
which implies $S=T$, $\P$-a.s., by the definition of $S_n$. Letting $n\to\infty$
in \eqref{stimasonesseenne} and \eqref{stimasonesseennebis}
 we conclude that
\eqref{stimaprlineare-dt} holds, so that $(Y,Z)\in \mathbb{M}$.
\qed


\begin{theorem} \label{Th:ex-un-dt}
Suppose that Hypothesis \ref{hyp:BSDE-dt} holds for some $(t,x)\in [0,T]\times K$.

Then there exists a unique pair $(Y,Z)$ in $\M^{t,x}$ which solves the
BSDE \eqref{BSDE-dt}.
\end{theorem}

\noindent{\bf Proof.} To simplify notation we drop the superscripts $t,x$
and we write the proof in the case $t=0$.
We use  a fixed point argument. Define the map $\Gamma: \M \rightarrow \M$
 as follows: for $(U,V) \in \M$,
 $(Y,Z)=\Gamma (U,V)$ is defined as the unique solution in $\M$ to the equation
\begin{equation}\label{mappabackward-ds}
  Y_s+\int_s^T \int_K Z_r(y) \; q(dr \;dy)=
   g(X_T)+\int_s^Tf (r,X_r, U_r,V_r)\; dr,\qquad s\in [0,T].
\end{equation}
From the assumptions on $f$ it follows that $\E \int_0^T |f( s,X_s,U_s,V_s)|^2 ds<\infty $,
so by Lemma \ref{BSDEcasobanale-dt} there exists a unique
$(Y,Z)\in\M$  satisfying
 \eqref{mappabackward-ds} and
 $\Gamma$ is a
well defined map.

We show that $\Gamma$ is a contraction if $\M$ is endowed with the equivalent norm
\begin{equation} \label{C-norm}||(Y,Z)||_{\M}^{2}
:= C|Y|^2_{\beta}+ ||Z||^2_{\beta},
\end{equation}
where
$$
|Y|^2_{\beta}:=
\E\int_0^T e^{\beta s}|Y_s|^2ds,
\qquad
||Z||^2_{\beta}:=
 \E \int_0^T  \int_K e^{\beta s}|Z_s(y)|^2\nu(s,X_{s},dy)\,ds,
 $$
for some constants $C>0$ and $\beta >0$ sufficiently large, that will be determined in the sequel.

Let $(U^1,V^1)$, $(U^2,V^2)$ be two elements of $\mathbb{M}$ and let
$(Y^1,Z^1)$, $(Y^2,Z^2)$ the associated solutions.
Lemma \ref{BSDEcasobanale-dt} applies to the difference $\overline{Y}=Y^1-Y^2$, $\overline{Z}=Z^1-Z^2$, $\overline{f}_s=f(s,X_s,U^1_s,V^1_s) -f(s,X_s,U^2_s,V^2_s)$ and
\eqref{identitaenergia-dt} yields, noting that $\overline{Y}_T=0$,
$$
\begin{array}{l}\dis
\E e^{\beta s}| \overline{Y}_s|^2 +\beta\,\E\int_s^T  e^{\beta r}|\overline{Y}_r|^2 dr
+ \E\int_s^T \int_K e^{\beta r}|\overline{Z}_r(y)|^2 \nu(r,X_{r},dy)\,dr\\
\dis\qquad
=
  2\E \int_s^T  e^{\beta r} \overline{Y}_{r} \,\overline{f}_r dr,
  \qquad s\in [0,T].
\end{array}
$$
From the Lipschitz conditions of $f$ and elementary inequalities it follows that
$$
\begin{array}{l}\dis
 \beta\,\E\int_0^T  e^{\beta s}|\overline{Y}_s|^2 ds
+ \E\int_0^T \int_K e^{\beta s}|\overline{Z}_s(y)|^2 \nu(s,X_{s},dy)\,ds
\\
\dis\qquad\le
  2L\E \int_0^T  e^{\beta s} |\overline{Y}_{s}| \,
 \left(  \int_K |\overline{V}_s(y)|^2 \nu(s,X_{s},dy)\right)^{1/2}
   ds
  + 2L'\E \int_0^T  e^{\beta s} |\overline{Y}_{s}| \,|\overline{U}_s|\, ds
\\
\dis\qquad
\le
\alpha
\E \int_0^T \int_K e^{\beta s}|\overline{V}_s(y)|^2 \nu(s,X_{s},dy)\,ds+
\frac{L^2}{\alpha} \E \int_0^T e^{\beta s}|\overline{Y}_s |^2 \, ds
\\
\dis\qquad
\quad+\gamma L' \E \int_0^T e^{\beta s}|\overline{Y}_s |^2 \, ds+
\frac{L'}{\gamma}  \E \int_0^T e^{\beta s}|\overline{U}_s |^2 \, ds
\end{array}
$$
for every $\alpha>0$, $\gamma>0$. This can be written
$$
\left( \beta- \frac{L^2}{\alpha} - \gamma L' \right)\,|\overline{Y}|^2_\beta
+ \|\overline{Z}\|_\beta^2 \le
\alpha
\|\overline{V}\|^2_\beta +
\frac{L'}{\gamma}  |\overline{U}  |^2_\beta.
$$
If we choose $\beta> L^2 +2L'$, it is possible to find $\alpha\in (0,1)$ such that
$$
\beta>\frac{L^2}{\alpha} +\frac{2L'}{\sqrt\alpha}.
$$
If $L'=0$ we see that $\Gamma$ is an $\alpha$-contraction on
$\M$ endowed with the norm (\ref{C-norm}) for $C= \beta- ({L^2}/{\alpha})$. If $L'>0$ we choose $\gamma=1/\sqrt{\alpha}$
and obtain
$$
\frac{L'}{\sqrt\alpha} \,|\overline{Y}|^2_\beta
+ \|\overline{Z}\|_\beta^2 \le
\alpha
\|\overline{V}\|^2_\beta +
 {L'}{\sqrt{\alpha}}  |\overline{U}  |^2_\beta=
 \alpha \,
 \left( \frac{L'}{\sqrt\alpha} \,|\overline{U}|^2_\beta
+ \|\overline{V}\|_\beta^2 \right),
$$
so that
$\Gamma$   is an $\alpha$-contraction on
$\M$ endowed with the norm \eqref{C-norm} for $C=( {L'}/\sqrt{\alpha} )$. In all cases
there exists a unique fixed point which is the required unique solution to the BSDE \eqref{BSDE-dt}.
\qed

 Next we prove some estimates on the solutions of the BSDE, which
show in particular the continuous dependence upon the data. Let us
consider two solutions $(Y^1,Z^1)$, $(Y^2,Z^2)\in \M^{t,x}$ to
the BSDE \eqref{BSDE-dt} associated with the drivers $f^1$ and $f^2$
and final data $g^1$ and $g^2$, respectively, which are
assumed to satisfy Hypothesis \ref{hyp:BSDE-dt}. Denote
$\overline{Y}=Y^1-Y^2$, $\overline{Z}=Z^1-Z^2$, $\overline{g}_T=
g^1(X^{t,x}_T)-g^2(X^{t,x}_T)$, $ \overline{f}_s= f^1(s,X^{t,x}_s,Y^2_s,Z^2_s(\cdot))-
f^2(s,X^{t,x}_s,Y^2_s,Z^2_s(\cdot)).
 $

\begin{proposition} \label{a priori estimates}
Suppose that Hypothesis \ref{hyp:BSDE-dt} holds for some $t,x$.
Let $(\overline{Y},\overline{Z})$ be the processes defined above.
Then the a priori estimates hold:
\begin{equation}\label{stima-diff-y}
\begin{array}{l}\dis
\sup_{s \in[t,T]} \E^{t,x}|\overline{Y}_s|^2 +
\E^{t,x} \int_t^T |\overline{Y}|_s^2 \,ds +
\E^{t,x} \int_t^T \int_K |\overline{Z}|_s^2\nu(s,X_s,dy) \,ds
\\\dis\qquad
\leq  C\left(\E^{t,x}|\overline{g}_T|^2
+ \,\E^{t,x} \int_t^T |\overline{f}_s|^2 \,ds \right),
\end{array}
\end{equation}
where $C$ is a constant depending only on $T,L,L'$.
\end{proposition}

\noindent{\bf Proof.} Again we drop the superscripts $t,x$.
Arguing as in the proof of \eqref{identitaenergia-dt} we obtain
$$\begin{array}{l}
\E |\overline{Y}_s|^2    +
\dis \E\int_s^T \int_K |\overline{Z}_r(y)|^2 \nu(r,X_r,dy) \, dr \\
\qquad = \dis\E |\overline{g}_T|^2 +2\E \int_s^T  \overline{Y}_{r} (f^1(r,X_r,Y_r^1,Z^1_r)-f^2(r,X_r,Y_r^2,Z^2_r) \, dr.
\end{array}
$$
By the
Lipschitz property of the driver $f^1$ we get
$$\begin{array}{lll}
\E |\overline{Y}_s|^2
& \leq &\dis
\E |\overline{g}_T|^2 +2\E \int_s^T  |\overline{Y}_{r}|(|f^1(r,X_r,Y_r^1,Z^1_r)-f^1(r,X_r,Y_r^2,Z^2_r)| +|\overline{f}_r|) \,dr
\\
& \leq &\dis \E |\overline{g}_T|^2 +2L'\E \int_s^T |\overline{Y}_{r}|^2 \,dr
+2L\E \int_s^T  |\overline{Y}_{r}|\left\{\int_K|\overline{Z}_r|^2 \nu(r,X_r,dy)\right\}^{\frac{1}{2}}\,dr
\\&&\dis
+2\E \int_s^T
 |\overline{Y}_{r}||\overline{f}_r| \,dr,\\
& \leq &  \dis\E |\overline{g}_T|^2 +C \E \int_s^T |\overline{Y}_{r}|^2 \,dr
+ \frac{1}{2}\E \int_s^T  \int_K |\overline{Z}_r|^2  \nu(r,X_r,dy)\,dr
\\&&\dis
+ \E \int_s^T |\overline{f}_r|^2 \,dr,
\end{array}
$$
for some constant $C$.
Hence we deduce
\begin{eqnarray}\E |\overline{Y}_s|^2  &+&  \frac{1}{2}\E\int_s^T \int_K |\overline{Z}_r|^2 \nu(r,X_r,dy) \, dr \nonumber \\
\label{pergronwall}
&\leq& \E |\overline{g}_T|^2 +\E \int_s^T  |\overline{f}_r|^2 \,dr + C\,
\E \int_s^T |\overline{Y}_{r}|^2 \,dr
\end{eqnarray}
and by Gronwall's lemma we get
\begin{eqnarray*} \E |\overline{Y}_s|^2  \leq e^{C(T-s)}\left(\E |\overline{g}_T|^2 +\E \int_t^T  |\overline{f}_r|^2 \,dr \right)
\end{eqnarray*}
and the conclusion follows from \eqref{pergronwall}.
\qed

From the a priori estimates we deduce the continuous dependence of the solution upon the data:
\begin{corollary}
Suppose that Hypothesis \ref{hyp:BSDE-dt} holds for some $t,x$.
Let $(Y,Z)$ the unique solution in $\M^{t,x}$ of the BSDE (\ref{BSDE-dt}). Then there exists a positive constant $C$, depending only on $T,L,L'$, such that
\begin{equation}
\E^{t,x}\int_t^T |Y_s|^2 \, ds  +  \E^{t,x}\int_t^T  \int_K |Z_s(y)|^2 \nu(s,X_s,dy)\, ds
 \leq C\E^{t,x}\left[ |g(X_T)|^2 +\int_t^T  |f(s,X_s,0,0)|^2 \, ds \right].
\end{equation}
\end{corollary}

\noindent{\bf Proof.}
The thesis follows from Proposition \ref{a priori estimates} setting $f^1=f$, $g^1=g$,
$f^2=0$ and $g^2=0$.
\qed


\section{Non linear variants of the Kolmogorov equation}\label{nlk}

Let us assume that $\nu$ is a
a transition measure on $K$ satisfying
\eqref{lambdalimitato}. $X$ denotes the Markov process constructed in
section
\ref{notations}, satisfying conditions 1-6
in subsection \ref{jumpmarkov}
 as well as \eqref{jumpkernel}.

In this section it is our purpose to present some nonlinear variants
of the classical backward Kolmogorov equation associated
to the Markov  process $X$ and to show that their solution can be
represented probabilistically by means of an appropriate BSDE
of the type considered above.

 Suppose that two functions $f,g$ are given,
satisfying the assumptions of Hypothesis \ref{hyp:BSDE-dt} for every $t\in [0,T]$,
$x\in K$. The equation
\begin{equation}\label{NLK}
  v(t,x)
  =
  g(x)
  +
  \int_t^T
\call_s  v(s,x) \,ds
+\int_t^T f\Big(s,x,v(s,x),v(s,\cdot)-v(s,x) \Big)\,ds,
\qquad t\in [0,T],\, x\in K,
\end{equation}
with unknown function $ v:[0,T]\times K\to \R$, will be called
the non linear Kolmogorov equation.
Equivalently, one requires that
for every $x\in K$ the map
$t\mapsto v(t,x)$ is absolutely continuous on $[0,T]$
and
\begin{equation}\label{NLKdifferenziale}
\left\{\begin{array}{l}\dis
\partial_t v(t,x)+
\call_t v(t,x)
+ f\Big(t,x,v(t,x), v(t,\cdot)-v(t,x) \Big)=0,
\\
v(T,x)=g(x),
\end{array}
\right.
\end{equation}
where the first equality is understood to hold almost everywhere on $[0,T]$,
the set of points where it may fail possibly depending on $x$.

The classical Kolmogorov equation corresponds to the case $f=0$.

Under appropriate boundedness assumptions we have the following immediate result:

\begin{lemma}\label{theoremnlkbdd}
Suppose that $f,g$
verify Hypothesis \ref{hyp:BSDE-dt} and, in addition,
\begin{equation}\label{ippernlk}
\sup_{t \in [0,T],\;x\in K} (|g(x)|+|f(t,x,0,0)|)<\infty.
\end{equation}
Then
 the nonlinear Kolmogorov equation has a unique solution in the class
 of measurable bounded functions.
\end{lemma}

\noindent{\bf Proof.} The result is essentially
known (see for instance \cite{B}, Chapter VII,  Theorem T3), so we only sketch the proof.
In the space of bounded measurable real functions on $[0,T]\times K $ endowed
with the supremum norm one can define
a map $\Gamma$ setting
 $v=\Gamma (u)$ where
 $$
  v(t,x)
  =
  g(x)
  +
  \int_t^T
\call_s  u(s,x) \,ds
+\int_t^T f\Big(s,x,u(s,x),u(s,\cdot)-u(s,x) \Big)\,ds,
\qquad t\in [0,T],\, x\in K.
$$
Using the boundedness condition \eqref{lambdalimitato} and the
Lipschitz character of $f$, by standard estimates one can prove that
$\Gamma$ has a unique fixed point, which is the required solution.
\qed

Now we plan to   remove the boundedness assumption
\eqref{ippernlk}. On the functions $f,g$ we will only impose the
 conditions required in
Hypothesis \ref{hyp:BSDE-dt}, for every $t\in [0,T]$,
$x\in K$.

\begin{definition}\label{defsoluznlk}
We say that a measurable function $ v:[0,T]\times K\to \R$
 is a solution of the non linear Kolmogorov equation \eqref{NLK}
if,  for every $t\in [0,T]$, $x\in K$,
 \begin{enumerate}
\item
$\dis\E^{t,x}\int_t^T\int_K |v(s,y)-v(s,X_{s})|^2\nu(s,X_s,dy)\,ds<\infty$;
 \item
$\dis \E^{t,x}\int_t^T  | v(s,X_{s})|^2 \,ds<\infty$;
\item  \eqref{NLK} is satisfied.
 \end{enumerate}
\end{definition}

\begin{remark} {\em
Condition 1 is equivalent to the fact that $v(s,y)-v(s,X_{s-})$
belongs to $\call^{2}(p^{t})$.
 Conditions  1 and 2 together
are equivalent to the fact that the pair $
\{(v(s,X_s),v(s,y)-  v(s,X_{s-});$ $s\in [t,T],y\in K\}
$
belongs to the space $\M^{t,x}$; in particular they hold true
 for every
measurable bounded function $ v$.
}
\end{remark}
\begin{remark} {\em
We need to verify that for a function $v$ satisfying the conditions 1 and 2 above the
equation \eqref{NLK} is well defined.

We first note that for every $x\in K$ we have, $\P^{0,x}$-a.s.,
$$
\int_0^T\int_K |v(s,y)-v(s,X_{s})|^2\nu(s,X^{0,x}_s,dy)\,ds
+\int_0^T  | v(s,X_{s})|^2 \,ds
<\infty.
$$
We recall that the law of the first
jump time $T_1$ is exponential with variable rate, according to
 \eqref{jumpkerneldue}. It follows
that the set $\{\omega\in\Omega\,:\,T_1(\omega)>T\}$
has positive $\P^{0,x}$ probability, and on this set we have
$X_{s}(\omega)=x$. Taking such an $\omega$ we conclude that
$$
\int_0^T\int_K |v(s,y)-v(s,x)|^2\nu(s,x,dy)\,ds
+\int_0^T  | v(s,x)|^2 \,ds
<\infty,\qquad x\in K.
$$
Since we are assuming $\sup_{t,x}\nu(t,x,K)<\infty$,
it follows from H\"older's inequality that
\begin{eqnarray*}
\int_0^T  | \call_s v(s,x)|  \,ds&\le&
\int_0^T\int_K |v(s,y)-v(s,x)|\nu(s,x,dy)\,ds
\\
&\le &c\left(
\int_0^T\int_K |v(s,y)-v(s,x)|^2\nu(s,x,dy)\,ds\right)^{1/2}<\infty
\end{eqnarray*}
for some constant $c$ and for all $x\in K$.

Similarly, from our assumption $ \E^{t,x} \int_0^T |f(s,X_s,0,0) |^2 ds <\infty$
 we deduce, arguing again on the jump time $T_1$, that
$$
\int_0^T |f(s,x,0,0) |^2\, ds
<\infty,\qquad x\in K,
$$
and from the Lipschitz conditions on $f$ we conclude that
\begin{eqnarray*}
\int_0^T  |  f (s,x,v(s,x),v(s,\cdot)-v(s,x)  ) |  \,ds&\le&
c_1\!\!\left( \int_0^T |f(s,x,0,0) |^2\, ds\right)^\frac12
\!\!\!\!+c_2\!\!\left(
\int_0^T |v(s,x ) |^2\, ds\right)^\frac12
\\
&  &+c_3\!\!\left(
\int_0^T\int_K |v(s,y)-v(s,x)|^2\nu(s,x,dy)\,ds\right)^\frac12
<\infty
\end{eqnarray*}
for some constants $c_i$ and for all $x\in K$.

We have thus verified that all the terms occurring in
equation \eqref{NLK} are well defined.
}
\end{remark}

In the following,
a basic role will be played by the  BSDEs: $\P^{t,x}$-a.s.
 \begin{equation}\label{bsdepernlk}
   Y_s^{t,x}+\int_s^T\int_KZ_r^{t,x}(y)\,q^{t}(dr\,dy) =
g(X_T) +\int_s^T f(r,X_r,Y_r^{t,x},Z_r^{t,x}(\cdot))\,dr,
\qquad s\in [t,T],
 \end{equation}
with unknown processes
$(Y_s^{t,x}, Z_s^{t,x})_{s\in [t,T]}$.  For every $(t,x)\in [0,T]\times K$
there exists a unique solution in the sense of theorem  \ref{Th:ex-un-dt}.
Note that $Y_t^{t,x}$ is deterministic.

We are ready to state the main result of this section.

\begin{theorem}\label{theoremnlk}
Suppose that hypothesis \ref{hyp:BSDE-dt} holds for every $(t,x)\in [0,T]\times K$.
Then the non linear Kolmogorov equation \eqref{NLK} has a unique solution $v$.

Moreover, for every
$(t,x)\in [0,T]\times K$ we have
\begin{eqnarray}
    Y_s^{t,x}&=&v(s,X_s),\label{identificazioneymarkov}
\\
Z_s^{t,x}(y)&= &v(s,y)-  v(s,X_{s-}),\label{identificazionezmarkov}
\end{eqnarray}
so that
in particular  $v(t,x)=Y_t^{t,x}$.
\end{theorem}

\begin{remark}\label{identificazionecommenti}\begin{em}
The  equalities \eqref{identificazioneymarkov} and
\eqref{identificazionezmarkov}
are understood as follows.
\begin{itemize}
\item
$\P^{t,x}$-a.s.,
equality \eqref{identificazioneymarkov} holds  for all $s\in [t,T]$.

Since the trajectories of $X$ are piecewise constant and cadlag
this is equivalent to the condition $\E^{t,x} \int_t^T|Y_s^{t,x}-v(s,X_s)|^2ds=0$.

\item
The equality
\eqref{identificazionezmarkov} holds for almost all $(\omega,s,y)$
with respect to the measure
\\
$ \nu(s,X_{s-}^{t,x}(\omega),dy)\,\P^{t,x}(d\omega)\,ds$,
i.e.
$$\E^{t,x} \int_t^T|Z_s^{t,x}(y)-v(s,y)+  v(s,X_{s-})|^2\nu(s,X_{s-},dy)\,ds=0.
$$
\end{itemize}
\end{em}
\end{remark}

\noindent{\bf Proof.} {\em Uniqueness.}
Let $v$ be a solution. It follows from equality \eqref{NLK} itself that
 $t\mapsto v(t,x)$ is absolutely continuous on $[0,T]$
for every $x\in K$. Since we assume that   the process
$v(s,y)-v(s,X_{s-})$
belongs to $\call^{2}(p^{t})$, we are in a position to apply the Ito formula
\eqref{Itoformula} to the
process $v(s,X_s)$, $s\in [t,T]$, obtaining, $\P^{t,x}$-a.s.,
$$
\begin{array}{lll}
v(s,X_s) &=&\dis
v(t,x)+\int_t^s\Big(\partial_r v(r,X_r) +\call_r v(r,X_r)\Big)\, dr
\\&&\dis +
\int_t^s \int_K \Big(v(r,y)-v(r,X_{r-})\Big)\, q^{t}(dr\,dy),
\qquad s\in [t,T].
\end{array}
 $$
 Taking into account that $v$ satisfies
\eqref{NLKdifferenziale} and that
$X$ has piecewise constant trajectories
 we obtain, $\P^{t,x}$-a.s.,
$$
\partial_r v(r,X_r)+
\call_r v(r,X_r)
+ f\Big(r,X_r,v(r,X_r), v(r,\cdot)-v(r,X_r) \Big)=0,
$$
for almost all $r\in [t,T]$.
It follows that, $\P^{t,x}$-a.s.,
\begin{equation}\label{perrelfond}
\begin{array}{lll}
v(s,X_s) &=&\dis
v(t,x)-\int_t^s f\Big(r, X_r, v(r,X_r),  v(r,\cdot)-v(r,X_r)\Big)\, dr
\\&&\dis +
\int_t^s \int_K \Big(v(r,y)-v(r,X_{r-})\Big)\, q^{t}(dr\,dy),
\qquad s\in [t,T].
\end{array}
\end{equation}
 Since $v(T,x)=g(x)$ for all $x\in K$, simple passages show that
$$
\begin{array}{l}
\dis v(s,X_s) + \int_s^T \int_K \Big(v(r,y)-v(r,X_{r-})\Big)\, q^{t}(dr\,dy)
\\\dis\qquad\qquad
=
g( X_T)
+
 \int_s^T f\Big(r, X_r, v(r,X_r),  v(r,\cdot)-v(r,X_r)\Big)\, dr,
\qquad s\in [t,T].
\end{array}
 $$
Therefore the pairs
$
 (Y_s^{t,x},Z_s^{t,x}(y))$ and $
(v(s,X_s),v(s,y)-  v(s,X_{s-})
$
are both solutions to the same BSDE under $\P^{t,x}$,
and therefore they coincide as members
of the space $\M^{t,x}$. The required equalities \eqref{identificazioneymarkov}
 and \eqref{identificazionezmarkov}
follow.
 In particular we have $v(t,x)=Y_t^{t,x}$, which proves uniqueness of the solution.

\medskip

{\em Existence.}
By theorem \ref{Th:ex-un-dt}, for every $(t,x)\in [0,T]\times K$ the  BSDE
 \eqref{bsdepernlk}
   has a unique solution
$(Y_s^{t,x}, Z_s^{t,x})_{s\in [t,T]}$ and, moreover,
$Y_t^{t,x}$ is deterministic, i.e. there exists a real number, denoted $v(t,x)$, such that
$\P^{t,x}(Y_t^{t,x}=v(t,x))=1$.

 We proceed by an approximation argument.
Let $f^n=(f\wedge n)\vee(-n)$, $g^n=(g\wedge n)\vee(-n)$ denote
the truncations of $f$ and $g$ at level $n$. By lemma \ref{theoremnlkbdd} there exists
a unique bounded measurable solution $v^n$ to the equation: for
$t\in [0,T]$,  $ x\in K$,
\begin{equation}\label{nlkn}
      v^n(t,x)
  =
  g^n(x)
  +
  \int_t^T
\call_s  v^n(s,x) \,ds
+\int_t^T f^n\Big(s,x,v^n(s,x),v^n(s,\cdot)-v^n(s,x) \Big)\,ds.
\end{equation}
By the first part of the proof, we known that
$$
v^n(t,x)=Y_t^{t,x,n},\qquad
v^n(s,X_s)=
Y_s^{t,x,n},\qquad
v^n(s,y)-  v^n(s,X_{s-})=
Z_s^{t,x,n}(y),
$$
in the sense of  remark \ref{identificazionecommenti}, where
$(Y_s^{t,x,n}, Z_s^{t,x,n})_{s\in [t,T]}$ is the unique solution to
 the  BSDE
 $$   Y_s^{t,x,n}+\int_s^T\int_KZ_r^{t,x,n}(y)\,q^{t}(dr\,dy) =
g^n(X_T) +\int_s^T f^n(r,X_r,Y_r^{t,x,n},Z_r^{t,x,n}(\cdot))\,dr,
\quad s\in [t,T].
$$
Comparing with \eqref{bsdepernlk} and applying Proposition \ref{a priori estimates} we deduce
that for some constant $c$
\begin{equation}\label{treterminiazero}
\begin{array}{l}\dis
 \sup_{s\in [t,T]}\E^{t,x} |Y_s^{t,x}-Y_s^{t,x,n}|^2
 +
 \E^{t,x}\int_t^T   |Y_s^{t,x}-Y_s^{t,x,n}|^2 ds
 \\\qquad\dis
+ \E^{t,x}\int_t^T \int_K  |Z_s^{t,x}-Z_s^{t,x,n}|^2 \nu(s,X_{s},dy)\,ds
\\
\dis\le
c \E^{t,x}|g(X_T)-g^n(X_T)|^2 \!+
c \E^{t,x}\int_t^T |   f (s,X_s,Y_s^{t,x },Z_s^{t,x }(\cdot))-
f^n(s,X_s,Y_s^{t,x },Z_s^{t,x }(\cdot))|^2\,ds\to 0
\end{array}
\end{equation}
where the right-hand side tends to zero by monotone convergence.

In particular it follows that
$$
|v(t,x)- v^n(t,x)|^2 = |Y_t^{t,x}-Y_t^{t,x,n}|^2
 \le
 \sup_{s\in [t,T]}\E |Y_s^{t,x}-Y_s^{t,x,n}|^2 \to 0,
 $$
which shows that $v$ is a measurable function. An application of the Fatou lemma
gives
$$\begin{array}{l}\dis
\E^{t,x} \int_t^T|Y_s^{t,x}-v(s,X_s)|^2ds+
\E^{t,x} \int_t^T|Z_s^{t,x}(y)-v(s,y)+  v(s,X_{s-})|^2\nu(s,X_{s-},dy)\,ds
\\\dis
\le \liminf_{n\to\infty}
\E^{t,x} \int_t^T|Y_s^{t,x}-v^n(s,X_s)|^2ds
\\\dis\qquad
+
\liminf_{n\to\infty}
\E^{t,x} \int_t^T|Z_s^{t,x}(y)-v^n(s,y)+  v^n(s,X_{s-})|^2\nu(s,X_{s-},dy)\,ds
\\\dis
=
\liminf_{n\to\infty}
\E^{t,x} \int_t^T|Y_s^{t,x}-Y_s^{n,t,x}|^2ds
+
\liminf_{n\to\infty}
\E^{t,x} \int_t^T|Z_s^{t,x}(y)-Z_s^{n,t,x}(y)|^2\nu(s,X_{s-},dy)\,ds=0
\end{array}
$$
by \eqref{treterminiazero}. This proves that
\eqref{identificazioneymarkov} and
 \eqref{identificazionezmarkov} hold.
  These formulae also imply that
$$
\begin{array}{l}\dis
 \E^{t,x}\int_t^T\int_K |v(s,y)-v(s,X_{s})|^2\nu(s,X_s,dy)\,ds+
  \E^{t,x}\int_t^T  | v(s,X_{s})|^2 \,ds
  \\\dis
  \qquad =
  \E^{t,x}\int_t^T\int_K |Z_s^{t,x}|^2\nu(s,X_s,dy)\,ds+
  \E^{t,x}\int_t^T  |Y_s^{t,x}|^2 \,ds  <\infty,
   \end{array}
 $$
 according to the requirements of
 definition \ref{defsoluznlk}. It only remains to show that $v$ satisfies
 \eqref{NLK}. This will follow from a passage to the limit in
 \eqref{nlkn}, provided we can show that
$$
  \int_t^T
\call_s  v^n(s,x) \,ds\to   \int_t^T
\call_s  v(s,x) \,ds,
$$
\begin{equation}\label{convergefn}
\int_t^T f^n\Big(s,x,v^n(s,x),v^n(s,\cdot)-v^n(s,x) \Big)\,ds\to
\int_t^T f\Big(s,x,v(s,x),v(s,\cdot)-v(s,x) \Big)\,ds.
\end{equation}

We first consider
 $$
 \begin{array}{l}\dis
 \E^{t,x}\left|
 \int_t^T
\call_s  v(s,X_s) \,ds- \int_t^T
\call_s  v^n(s,X_s) \,ds \right|
\\\dis\qquad
=  \E^{t,x}\left|
 \int_t^T\int_K[v(s,y)-v(s,X_s )
 - v^n(s,y)+v^n(s,X_s ) ]\,\nu(s, X_s ,dy)
 \,ds \right|
 \\\dis\qquad
 =  \E^{t,x}\left|
 \int_t^T\int_K(Z_s^{t,x}-Z_s^{t,x,n})\,\nu(s, X_s ,dy)
 \,ds \right|
\\\dis\qquad
\le (T-t)^{1/2} \sup_{t,x}\nu(t,x,K)^{1/2}
\left( \E^{t,x}\int_t^T \int_K  |Z_s^{t,x}-Z_s^{t,x,n}|^2 \nu(s,X_{s},dy)\,ds
\right)^{1/2}
\end{array}
$$
which tends to zero, by \eqref{treterminiazero}. So for a subsequence (still denoted $v^n$)
we have
$ \int_t^T
\call  v^n(s,X_s) \,ds$ $\to$ $\int_t^T
\call  v(s,X_s) \,ds $ $\P^{t,x}$-a.s.
Note that, according to \eqref{jumpkerneldue},  the first
jump time $T_1^{t}$ has exponential law, so the set
$\{\omega\in\Omega\,:\,T_1^{t}(\omega)>T\}$
has positive $\P^{t,x}$ probability,   and on this set we have
$X_{s}(\omega)=x$.
Taking such an $\omega$ we conclude that
$ \int_t^T
\call_s  v^n(s,x) \,ds \to \int_t^T
\call_s  v(s,x) \,ds $.

 To prove \eqref{convergefn} we compute
 $$
 \begin{array}{l}\dis
 \E^{t,x}\bigg| \int_t^T f\Big(s,X_s,v(s,X_s),v(s,\cdot)-v(s,X_s) \Big)\,ds
 \\\dis\qquad\qquad
 -
 \int_t^T f^n\Big(s,X_s,v^n(s,X_s),v^n(s,\cdot)-v^n(s,X_s) \Big)\,ds \bigg|
  \\\dis
  =  \E^{t,x}\bigg| \int_t^T [f (s,X_s,Y_s^{t,x},Z_s^{t,x}) -
 f^n (s,X_s,Y_s^{t,x,n},Z_s^{t,x,n}  )]\,ds \bigg|
 \\\dis
  \le  \E^{t,x}  \int_t^T |f (s,X_s,Y_s^{t,x},Z_s^{t,x}) -
  f^n (s,X_s,Y_s^{t,x},Z_s^{t,x}  )|\,ds
   \\\dis\qquad\qquad
  +
  \E^{t,x}  \int_t^T |f^n (s,X_s,Y_s^{t,x},Z_s^{t,x}) -
 f^n (s,X_s,Y_s^{t,x,n},Z_s^{t,x,n}  )|\,ds  .
 \end{array}
 $$
 In the right-hand side, the first integral tends to zero by monotone convergence.
 Since $f^n$ is a truncation of $f$, it satisfies the Lipschitz condition \eqref{flipschitz} with
 the same constants $L,L'$ independent of $n$; therefore the second integral can be estimated
 by
 $$\begin{array}{l}\dis
 L'\, \E^{t,x}  \int_t^T | Y_s^{t,x}-Y_s^{t,x,n} |\,ds  + L\,
 \E^{t,x}  \int_t^T \left(\int_K| Z_s^{t,x}(y)- Z_s^{t,x,n} (y )|^2 \,\nu(s,X_s,dy)
 \right)^{1/2}\,ds
 \\\dis
 \le
 L'\, \left( (T-t)\,\E^{t,x}  \int_t^T | Y_s^{t,x}-Y_s^{t,x,n} |^2\,ds \right)^{1/2}
 \\\dis \qquad+
  L\,
\left( (T-t)\,
 \E^{t,x}  \int_t^T  \int_K| Z_s^{t,x}(y)- Z_s^{t,x,n} (y )|^2 \,\nu(s,X_s,dy)\,ds
 \right)^{1/2},
  \end{array}
   $$
   which tends to zero, again by \eqref{treterminiazero}.
    So for a subsequence (still denoted $v^n$)
we have
$$
  \int_t^T f^n\Big(s,X_s,v^n(s,X_s),v^n(s,\cdot)-v^n(s,X_s) \Big)\,ds
 \to
   \int_t^T f\Big(s,X_s,v(s,X_s),v(s,\cdot)-v(s,X_s) \Big)\,ds
   $$
   $\P^{t,x}$-a.s.
Picking an $\omega$ in the set
$\{\omega\in\Omega\,:\,T_1^{t}(\omega)>T\}$ as before
we conclude that \eqref{convergefn} holds, and the proof is finished.
\qed


\section{Optimal control}\label{sec-control}

\subsection{Formulation of the problem}

In this section we start again with
a measurable space $(K,\calk)$  and
a transition measure $\nu$  on $K$, satisfying
\eqref{lambdalimitato}. The process $X$ is constructed as described
in section \ref{notations}.

The data specifying the optimal control problem that
we will address are a measurable space
$(U,\calu)$, called the
  action (or decision) space,
 a running cost function $l$, a (deterministic, finite)
 time horizon $T>0$, a terminal cost function $g$,
and another function $r$ specifying the effect of the control
process.

For every $t\in[0,T]$
we define an admissible  control process, or simply a control, as an
$\F^{t}$-predictable process $(u_s)_{s\in[t,T]}$ with values in $U$.
The set of admissible
control processes is denoted $\cala^{t}$.

We will make the following assumptions.

\begin{hypothesis}\label{hyp:controllo}
\begin{enumerate}
\item  $(U,\calu)$ is a measurable space.
\item $r:[0,T]\times K\times K\times U\to \R$
is $\calb([0,T])\otimes \calk \otimes \calk\otimes \calu$-measurable and
 there exist a
constant $C_r> 0$  such that
\begin{equation}\label{ellelimitato}
0\le r(t,x,y,u)\le C_r,\qquad t\in [0,T],
\,x,y \in K, u\in U.
\end{equation}

\item
$g:K\to\R$ is $\calk$-measurable and
\begin{equation}\label{gnonlimitato}
 \E^{t,x}|g(X_T)|^2<\infty,
\qquad t \in [0,T],\;x\in K.
\end{equation}
\item
$l:[0,T]\times K\times U\to \R$
is $\calb([0,T])\otimes \calk\otimes \calu$-measurable, and there
exists $\alpha>1$ such that for every
$t\in [0,T]$,
$x\in K$ and $u(\cdot)\in \cala^{t} $ we have
\begin{equation}\label{lnonlimitato}
\inf_{u\in U}l(t,x,u)>-\infty,\qquad
\E^{t,x}\int_t^T |\inf_{u\in U}l(s,X_s,u)|^2\,ds<\infty,
\end{equation}
\begin{equation}\label{elleammiss}
\E^{t,x}\left(\int_t^T |l(s,X_s,u_s)|\,ds\right)^\alpha<\infty.
\end{equation}

\end{enumerate}
\end{hypothesis}

\begin{remark}\begin{em}
We note that the cost functions $g$ and $l$ need not be bounded.
Clearly, \eqref{elleammiss} follows from the other assumptions
if we assume for instance that
$\E^{t,x}\int_t^T |\sup_{u\in U}l(s,X_s,u)|\,ds<\infty$ for all
$t\in [0,T]$  and
$x\in K$.
\end{em}
\end{remark}

To any $(t,x)\in [0,T]\times K$
  and any control
  $u(\cdot)\in\cala^{t}$ we associate a probability measure
  $\P_u^{t,x}$ on $(\Omega,\calf)$ by a change of measure
  of Girsanov type, as we now describe. Recalling the definition of the
  jump times $T_n^{t}$ in \eqref{deftn}
we  define,
for $ s\in [t,T]$,
$$
L_s^{t}=
\exp\left(\int_t^s\int_K (1-r(z,X_z,y,u_z))\,\nu(z,X_z,dy)\,dz\right)
\prod_{n\ge1\,:\,T_n^{t}\le s}r(T_n^{t},X_{T_{n-}^{t}},X_{T_n^{t}},u_{T_n^{t}}),
$$
with the convention that the last product equals $1$ if there are
no indices $n\ge 1$ satisfying $T_n^{t}\le s$. It is a well-known
result that $L^{t}$ is a nonnegative supermartingale relative to
$\P^{t,x}$ and $\F^{t}$  (see \cite{J}
Proposition 4.3, or \cite{Bo-Va-Wo}), solution to the equation
$$
L_s^{t}=
1+\int_t^s\int_K L_{z-}^{t}\,(r(z,X_{z-},y,u_z)-1)\;q^t(dz\,dy),
\qquad s\in [t,T].
$$
As a consequence of the boundedness assumption in
(\ref{lambdalimitato}) it
can be proved, using for instance Lemma 4.2 in \cite{Con-Fuh},
or \cite{B} Chapter VIII Theorem T11, that
for every $\gamma>1$ we have
\begin{equation}\label{Lintegrabile}
 \E^{t,x} [|L_T^{t}|^\gamma]<\infty,
\qquad \E^{t,x} L_T^{t}=1,
\end{equation}
 and therefore
 the process
$L^{t}$ is a martingale (relative to
$\P^{t,x}$ and $\F^{t}$).
Defining a probability $\P_u^{t,x}$ by
$\P_u^{t,x}(d\omega)=L_T^{t}(\omega)\P^{t,x}(d\omega)$,
we introduce the cost functional corresponding to $u(\cdot)\in\cala^{t}$ as
$$
J(t,x,u(\cdot))=\E_u^{t,x}\left[
\int_t^Tl(s,X_s,u_s)\,ds + g(X_T)\right],
$$
where
$\E_u^{t,x}$ denotes the
expectation under $\P_u^{t,x}$. Taking into account \eqref{gnonlimitato},
\eqref{elleammiss}, \eqref{Lintegrabile} and using the H\"older inequality it is easily
seen that  the cost is finite for every
admissible control. The control problem starting at
$(t,x)$ consists in minimizing $
J(t,x,\cdot)$ over  $\cala^{t}$.

We finally introduce the
value function
$$
v(t,x)= \inf_{u(\cdot)\in\cala^{t} }J(t,x,u(\cdot)),
\qquad t\in [0,T], \,x\in K.
$$

The previous formulation of the optimal control problem by means of a change
of probability measure is classical (see e.g.  \cite{ElK},  \cite{E},  \cite{B}).
Some comments may be useful
at this point.
\begin{remark}\begin{em}
\begin{enumerate}
\item
We recall (see e.g. \cite{B}, Appendix A2, Theorem T34) that a
process $u$ is $\F^{t}$-predictable if and only if it admits
the representation
\begin{equation}\label{rappresentazione-predictable-proc}
u_s(\omega )= \sum_{n \geq 0} u^{(n)}_s(\omega )\,1_{T_n^{t}(\omega) <s
\leq T_{n+1}^{t}(\omega)}
\end{equation}
where for each $n \ge 0$ the mapping $(\omega,s) \mapsto
u^{(n)}_s(\omega )$ is $\calf_{[t,T_n^{t}]} \otimes \calb
([t,\infty))$-measurable. Moreover,
$\calf_{[t,T_n^{t}]}=\sigma(T_i^{t},
X_{T_i^{t}}, 0\le i\le n)$ (see e.g. \cite{B},
Appendix A2, Theorem T30). Thus  the fact that   controls are predictable processes
admits the following  interpretation:   at each
time  $T_n^{t}$ (i.e., immediately after a jump) the controller,
having observed  the random variables
$T_i^{t},X_{T_i^{t}}$  $(0\le i\le n)$,  chooses his current control
action, and updates her/his decisions only at time $T^{t}_{n+1}$.

\item
 It can be proved (see
\cite{J} Theorem 4.5) that the compensator  of $p^{t}(ds\,dy)$
under $\P_u^{t,x}$ is
$$
r(s,X_{s-},y,u_s)\, \nu(s,X_{s-},dy)\,ds,
$$
whereas the  compensator of $p^{t}(ds\,dy)$
under $\P^{t,x}$ was $\nu(s,X_{s-},dy)\,ds$.
This explains that  the choice of a given control $u(\cdot)$ affects the
stochastic system multiplying
its  compensator  by $r(s,x,y,u_s)$.

\item
  We call control law an arbitrary measurable function
$\underline{u} : [0,T]\times K \to U$. Given a control law
one can define an admissible control $u$ setting
$
u_s=\underline{u}(s,X_{s-}).
$ Controls of this form are called feedback controls. For a feedback control
the compensator  of $p^{t}(ds\,dy)$ is
$
r(s,X_{s-},y,\underline{u}(s,X_{s-}))\, \nu(s,X_{s-},dy)\,ds
$ under $\P_u^{t,x}$.
Thus, in this case the controlled system  is
 a Markov process corresponding to the
transition measure
\begin{equation}\label{rateunderlaw}
    r(s,x,y,\underline{u}(s,x))\, \nu(s,x,dy)
\end{equation}
instead of
$\nu(s,x,dy)$.

We will see later that
an optimal control can often be found  in feedback form. In this case,
 even if the original
process was time-homogeneous (i.e. $\nu$ did not depend on time)
the optimal process is not, in general, since the control law
may depend on time.
\end{enumerate}
\end{em}
\end{remark}

\begin{remark}\label{confrontocontrollo}\begin{em}
Our formulation of the optimal control problem should be compared
with another classical approach (see e.g. \cite{Guo-He}, \cite{Da-bo})
that we describe informally.
One may start  with the same running and terminal cost functions
$l,g$ as before, but with
 a  jump rate function $\lambda^u(t,x)$ and a
jump measure $\pi^u(t,x,A)$ which also depend on the control parameter $u\in U$
as well as on
$t\in[0,T]$, $x\in K$, $A\in \calk$.
Controls only consist in  feedback laws, i.e. functions
$
\underline{u} : [0,T]\times K \to U.
$
Given any such $\underline{u}(\cdot,\cdot) $ one constructs a jump Markov process,
on some probability space,
with  jump rate function and jump measure given, respectively,  by
$$\lambda^{\underline{u}(t,x)}(t,x),
\qquad
\pi^{\underline{u}(t,x)}(t,x,A),
$$
or, equivalently, with rate measure $\lambda^{\underline{u}(t,x)}(t,x)\,
\pi^{\underline{u}(t,x)}(t,x,A)$.
Thus, together with the initial state and starting time, the choice of a
control law $\underline{u}(\cdot,\cdot)$ determines
the law of the process and consequently the corresponding cost
that we now denote
$J(\underline{u})$ (the cost functional
being defined in terms of $l$ and $g$ similarly as before).

Under appropriate conditions this optimal control problem can be reduced to our setting.
For instance suppose that there exist (fixed) jump rate function and jump measure
$\lambda(t,x), \pi(t,x,A)$ (equivalently, a rate measure
$\nu(t,x,A)=\lambda(t,x)\, \pi(t,x,A)$)
as in section \ref{notations} and that we have the implications
\begin{equation}\label{riduzcontrollo}
    \pi(t,x,A)=0\;\Rightarrow \;\pi^u(t,x,A)=0, \qquad
\lambda(t,x)=0\;\Rightarrow\; \lambda^u(t,x)=0,
\end{equation}
for every $t,x,A,u$.  Then denoting $y\mapsto r_0(x,t,y,u)$ the Radon-Nikodym derivative
of $\pi^u(t,x,\cdot)$ with respect to $\pi(t,x,\cdot)$, whose existence is granted by
\eqref{riduzcontrollo}, we can define
$$
r(t,x,y,u)=r_0(x,t,y,u)\,\frac{\lambda^u(t,x)}{\lambda(t,x)},
$$
with the convention that $0/0=1$.
Suppose also  that $r$ is measurable and bounded, so that it satisfies
Hypothesis \ref{hyp:controllo}-2.  Then we have the identity
$$
r(t,x,y,u)\, \nu(t,x,A)=
r(t,x,y,u)\, \lambda(t,x)\, \pi(t,x,A)=
\lambda^{u}(t,x)\,
\pi^{u}(t,x,A),
$$
whence it follows that the choice of any
control law $\underline{u}(\cdot,\cdot)$, giving rise to the rate
measure \eqref{rateunderlaw}, will correspond to a cost equal to
$J(\underline{u})$. Therefore the required reduction is in fact possible.

We mention however that, unless some condition like
\eqref{riduzcontrollo} is verified,
 the class of control problems specified
by the initial data $\lambda^u(t,x)$ and $\pi^u(t,x,A)$ is in general
larger than the one we address in this paper. This can be seen noting
that in our framework all the controlled processes have laws which
are absolutely continuous with respect to a single uncontrolled process
(the one corresponding to $r\equiv 1$) whereas this might not be the
case for the rate measures $\lambda^{\underline{u}(t,x)}(t,x)\,
\pi^{\underline{u}(t,x)}(t,x,A)$ when $\underline{u}(\cdot,\cdot)$ ranges in
the set of all possible control laws: a precise verification might be based
on the results of Section 4 in \cite{J} where absolute continuity
of the laws of marked point processes is characterized in terms
of their compensators.
\end{em}
\end{remark}

\subsection{The Hamilton-Jacobi-Bellman equation and the solution to the control problem}

The Hamilton-Jacobi-Bellman  (HJB) equation is the following non linear Kolmogorov equation:
 for every $t\in [0,T],\, x\in K$,
\begin{equation}\label{HJBdet}
  v(t,x)
  =
  g(x)
  +
  \int_t^T
\call_s  v(s,x) \,ds
+\int_t^T f\Big(s,x,v(s,\cdot)-v(s,x) \Big)\,ds,
\end{equation}
where $\call_s$ denotes
the generator of the Markov process $X$ as before, and
$f$ is the  hamiltonian function defined by
\begin{equation}\label{defhamiltonian}
    f( s,x,z(\cdot))=\inf_{u\in U}\left\{
l(s , x,u)+ \int_K z(y) \, (r
(s,x,y,u)-1)\,\nu (s,x,dy)\right\},
\end{equation}
for  $s \in [0,T]$, $x\in K$,
 $z\in L^2(K,\mathcal{K}, \nu(s,x,dy) )$. The (possibly empty) set of minimizers
 will be denoted
 \begin{equation}\label{minimizzanti}
    \Gamma(s,x, z(\cdot))=
    \{u\in U\,:\,
     f( s,x,z(\cdot))=
l(s , x,u)+ \int_K z(y) \, (r
(s,x,y,u)-1)\,\nu (s,x,dy)\}.
\end{equation}
We note that the HJB equation can be written in the alternative form:
$$
\left\{\begin{array}{l}\dis
\partial_t v(t,x)+
\inf_{u\in U}\left\{
l(t, x,u)+ \int_K (v(t,y)-v(t,x)) \, r
(t,x,y,u) \,\nu (t,x,dy)\right\} =0,
\\
v(T,x)=g(x), \qquad t\in [0,T],\, x\in K,
\end{array}
\right.
$$
but we will rather use \eqref{HJBdet} in order to make a connection with previous
results.
To study the HJB equation
we use the notion of solution presented in definition \ref{defsoluznlk}. We have the following preliminary result:

\begin{lemma}\label{ipotesiverificate}
Under Hypothesis \ref{hyp:controllo} the assumptions of Hypothesis
\ref{hyp:BSDE-dt}  hold true for every $(t,x)\in [0,T]\times K$ and consequently
the HJB equation has a unique solution according to
 Theorem \ref{theoremnlk}.
\end{lemma}

\noindent{\bf Proof.}  Hypothesis \ref{hyp:BSDE-dt}-1 and 4 coincide with
 \eqref {gnonlimitato} and \eqref{lnonlimitato} respectively. The only non trivial
 verification is the Lipschitz condition \eqref{flipschitz}: this follows from the boundedness
 assumption
  \eqref{ellelimitato} which implies that, for every $s \in [0,T]$, $x\in K$,
 $z,z' \in L^2(K,\mathcal{K}, \nu(s,x,dy) )$, $u\in U$,
  $$ \begin{array}{l}
\dis\int_K z(y) \, (r
(s,x,y,u)-1)\,\nu (s,x,dy) \\\dis\;
\le
  \int_K |z(y)-z'(y)| \, (r
(s,x,y,u)-1)\,\nu (s,x,dy) + \int_K z'(y) \, (r
(s,x,y,u)-1)\,\nu (s,x,dy)
\\\dis\;
\le
 (C_r+1)\,\nu (s,x,K)^{1/2}\left( \int_K |z(y)-z'(y)|^2 \,\nu (s,x,dy)\right)^{1/2} + \int_K z'(y) \, (r
(s,x,y,u)-1)\,\nu (s,x,dy) ,
\end{array}
  $$
  so that adding $l(s , x,u)$ to both sides and taking the infimum over $u\in U$
  it follows that
  $$
  f( s,x,z(\cdot))\le L \left( \int_K |z(y)-z'(y)|^2 \,\nu (s,x,dy)\right)^{1/2} +
  f( s,x,z'(\cdot))
  $$
  where  $L= (C_r+1)\,\sup_{t,x}\nu (t,x,K)^{1/2}<\infty$; exchanging $z$ and $z'$
  we obtain \eqref{flipschitz}.
  \qed

We can now state our main result.

\begin{theorem}\label{hjbrisolvecontrollo} Suppose that
 Hypothesis \ref{hyp:controllo}  holds.

Then there exists a unique solution $v$  to the
HJB equation.
Moreover, for any $t\in [0,T]$, $x\in K$ and
any admissible control $u(\cdot)\in \cala^{t}$ we have
$v(t,x)\le J(t,x,u(\cdot))$.

Suppose in addition that the sets $\Gamma$ introduced in
\eqref{minimizzanti}
are non empty and for every $t\in [0,T]$, $x\in K$ one can find
 an $\F^t$-predictable process $u^{*,t,x}(\cdot)$ in $U$
satisfying
\begin{equation}\label{minselector}
    u^{*,t,x}_s \in \Gamma(s,X_{s-},v(s,\cdot)-  v(s,X_{s-})  ),
\end{equation}
$\P^{t,x}$-a.s. for almost all $s\in [t,T]$.

Then
$u^{*,t,x}(\cdot)\in \cala^{t}$, it is an optimal control, and $v(t,x)$ coincides
with the value function, i.e.
$ v(t,x)=J(t,x,u^{*,t,x}(\cdot))$.
\end{theorem}

\begin{remark}\begin{em}
\begin{enumerate}
  \item The existence of
a process $u^{*,t,x}$ satisfying \eqref{minselector} is crucial in order
to apply the theorem and solve the optimal control problem
in a satisfactory way.  It is possible to formulate general sufficient
conditions for the existence of $u^{*,t,x}$: see Proposition
\ref{selezione} below. The proof of this proposition makes it clear that
in general the process $u^{*,t,x}$ may depend on $t,x$.
  \item
  Suppose that
there exists a measurable function
$\underline{u} :[0,T]\times K  \to U$ such that
\begin{equation}\begin{array}{l}\dis
l(s , x,\underline{u}( s,x))+ \int_K  \Big(v(s, y) -v(s, x)\Big)
 \, \Big(r (s,x,y,\underline{u}( s,x))-1\Big)\,\nu (s,x,dy)
 \\\dis\qquad
 =\inf_{u\in U}\left\{
l(s , x,u)+ \int_K  \Big(v(s, y) -v(s, x)\Big)  \, \Big(r
(s,x,y,u)-1\Big)\,\nu (s,x,dy)\right\},
\end{array}
\end{equation}
for all $s\in[0,T]$, $x\in K$, where $v$ denotes the solution of the HJB equation.
We note that in specific
situations it is possible to compute explicitly the function
$\underline{u}$.
 Then the process
$$
u^{*,t,x}_s= \underline{u}(s,X_{s-})
  $$
satisfies \eqref{minselector} and is therefore optimal.
Note that in this case
the optimal control is  in feedback form and the feedback law
$\underline{u}$
 is the same for every starting point $(t,x)$.

\end{enumerate}

\end{em}
\end{remark}

\noindent{\bf Proof.} Existence and uniqueness of a solution to the HJB equation,
in the sense of definition \ref{defsoluznlk}, is a consequence of lemma \ref{ipotesiverificate}
and Theorem \ref{theoremnlk}.
All the other statements of the theorem are immediately deduced
from the following identity, sometimes
called the fundamental relation:
for any $t\in [0,T]$, $x\in K$ and
any admissible control $u(\cdot)\in \cala^{t}$,
\begin{equation}\label{fundrel}
\begin{array}{lll}\dis
    v(t,x)&=&\dis J(t,x,u(\cdot)) +\E_u^{t,x} \int_t^T \bigg\{
     f \Big( s,X_s,v(s,\cdot)- v(s,X_s) \Big)
     \\&&\dis
     -
l(s , X_s,u_s)
 - \int_K(v(s,y)- v(s,X_s))\, (r
(s,X_s,y,u_s)-1)\,\nu ( s,X_s,dy)\bigg\}\,ds.
\end{array}
\end{equation}
Indeed, the term in curly brackets $\{\ldots\}$ is non positive by the definition
of the hamiltionian function \eqref{defhamiltonian}, and it equals zero when
$u(\cdot)$ coincides with $u^{*,t,x}(\cdot)$ by \eqref{minimizzanti}.

To finish the proof we   show that  \eqref{fundrel} holds. Applying
the Ito formula
\eqref{Itoformula} to the
process $v(s,X_s)$, $s\in [t,T]$,
and proceeding as in the proof of Theorem \ref{theoremnlk} we arrive at
equality
\eqref{perrelfond}, that we write for $s=T$:   recalling that
$v(T,x)=g(x)$ for all $x\in K$ we obtain
$$
\begin{array}{lll}
v(t,x) &=&\dis g(X_T)
+\int_t^T f\Big(s, X_s, v(s,X_s),  v(s,\cdot)-v(s,X_s)\Big)\, ds
\\&&\dis -
\int_t^T \int_K \Big(v(s,y)-v(s,X_{s-})\Big)\, q^{t}(ds\,dy).
\end{array}
$$
Since
$q^{t}(ds \, dy)= p^{t}(ds\, dy) - \nu(s,X_{s-},dy)\,ds$, we have,
adding and subtracting some terms,
$$
\begin{array}{lll}
v(t,x) &=&\dis g(X_T) +
\int_t^T
l(s , X_s,u_s) \,ds
\\&&\dis
+\int_t^T\bigg\{  f\Big(s, X_s, v(s,X_s),  v(s,\cdot)-v(s,X_s)\Big)
-l(s , X_s,u_s)
\\&&\dis\qquad
  - \int_K(v(s,y)- v(s,X_s))\, (r
(s,X_s,y,u_s)-1)\,\nu (s, X_s,dy)\bigg\}\,ds
\\&&\dis
- \int_t^T \int_K \Big(v(s,y)-v(s,X_{s-})\Big)\,
\Big(p^{t}(ds\, dy) - r(s,X_s,y,u_s)\nu(s,X_{s-},dy)\,ds\Big).
\end{array}
$$
Now \eqref{fundrel}  follows by taking the expectation with respect to $\P^{t,x}_u$,
provided we can show that the last term (the stochastic integral) has
mean zero with respect to $\P^{t,x}_u$. Since
the $\P^{t,x}_u$-compensator of $p^{t}(ds\, dy)$
is $r(s,X_{s-},y,u_s)\nu(s,X_{s-},dy)\,ds$, it is enough
to verify that the integrand $v(s,y)-v(s,X_{s-})$ belongs to $\call^1(p^{t})$
(with respect to $\P^{t,x}_u$), i.e. that the following integral, denoted $I$, is
finite:
$$
I= \E^{t,x}_u \int_t^T \int_K |v(s,y)-v(s,X_{s-})|\,
 r(s,X_s,y,u_s)\nu(s,X_{s},dy)\,ds.
 $$
 We have, by \eqref{ellelimitato} and the H\"older inequality,
$$
\begin{array}{lll}
  I & \le & \dis C_r\,
  \E^{t,x}_u   \int_t^T \int_K |v(s,y)-v(s,X_{s })|\,
  \nu(s,X_{s},dy)\,ds
   \\
    &  = & \dis C_r\,
  \E^{t,x}  \left[ L_T^{t} \int_t^T \int_K |v(s,y)-v(s,X_{s })|\,
  \nu(s,X_{s},dy)\,ds\right]
   \\
    &  \le  & \dis C_r \left(\E^{t,x}   | L_T^{t} |^2 \right)^\frac12
   \left( \E^{t,x} \left| \int_t^T \int_K |v(s,y)-v(s,X_{s })|\,
  \nu(s,X_{s},dy)\,ds\right|^2\right)^\frac12
  \\
    &  \le  & \dis C_r \left(\E^{t,x}   | L_T^{t} |^2 \right)^\frac12
   \left((T-t)  \sup_{t,x}\nu(t,x,K) \E^{t,x}\!\!   \int_t^T \!\int_K |v(s,y)-v(s,X_{s })|^2\,
  \nu(s,X_{s},dy)\,ds  \right)^\frac12\!\!.
\end{array}
$$
Recalling
\eqref{Lintegrabile} and the integrability condition in definition \ref{defsoluznlk}
we conclude that $I<\infty$ and this finishes the proof.
\qed

 As a consequence of theorem \ref{theoremnlk} we can also conclude that
 the value function and the optimal control law can also be  represented
 by means of the solution $(Y_s^{t,x}, Z_s^{t,x})_{s\in [t,T]}$ of the following
 BSDE:  $\P^{t,x}$-a.s.
   $$
    Y_s^{t,x}+\int_s^T\int_KZ_r^{t,x}(y)\,q^{t}(dr\,dy) =
g(X_T) +\int_s^T f(r,X_r,Z_r^{t,x}(\cdot))\,dr,
\qquad s\in [t,T],
$$
with fixed $(t,x)\in [0,T]\times K$ and the generator equal to
the hamiltonian function $f$. As before,
equalities  \eqref{identificazioneyezmarkov} below are understood as explained in
remark \ref{identificazionecommenti}.

\begin{corollary} \label{identifteorema}
Under the
assumptions of theorem  \ref{hjbrisolvecontrollo},
for every
$(t,x)\in [0,T]\times K$ we have
\begin{equation}\label{identificazioneyezmarkov}
    Y_s^{t,x}=v(s,X_s),
\quad
Z_s^{t,x}(y)= v(s,y)-  v(s,X_{s-}).
\end{equation}

In particular, the value function and an optimal control are given by the formulae
$$ v(t,x)=Y_t^{t,x},\qquad
u^{*,t,x}_s=\underline{u}(s,X_{s-},Z_{s}^{t,x}(\cdot)  ).
$$

\end{corollary}

As mentioned before,
general conditions can be
formulated for the existence of
a process $u^{*,t,x}$ satisfying \eqref{minselector}, hence of
an optimal control. This is done in  the
following proposition,  by means of an appropriate selection theorem.

\begin{proposition} \label{selezione}

In addition to the assumptions in Hypothesis \ref{hyp:controllo},
suppose that $U$ is a compact metric space  with its
Borel $\sigma$-algebra $\calu$ and that the functions
$r(s,x,\cdot), l (s,x,\cdot):U\to \R$ are
continuous for every $s\in [0,T]$, $x\in
K$. Then a process $u^{*,t,x}$ satisfying \eqref{minselector} exists
and all the conclusions of Theorem \ref{hjbrisolvecontrollo} hold true.
\end{proposition}

\noindent {\bf Proof.} We fix $t,x$ and consider
 the measure
$\mu(d\omega\,ds)= \P^{t,x}(d\omega)\, ds$ on the product
$\sigma$-algebra $\calg:=\calp^t\otimes \calb([t,T])$.
 Let $\bar\calg$ denote its $\mu$-completion and
consider the complete measure space $(\Omega\times [t,T],  \bar\calg, \mu)$.
Let $v$ denote the solution of the HJB equation.
 Define a map $F:\Omega\times [0,T] \times U\to\R $ setting
$$
F(\omega, s,u)=  l (s,  X_{s-}(\omega),u)
+ \int_K  \Big(v(s, y) -v(s, X_{s-}(\omega))\Big) \,
\Big(r(s,X_{s-},y,u)-1\Big)\,\nu(s, X_{s-}(\omega),dy).
$$
Then $ F(\cdot,\cdot, u)$ is $\bar\calg$-measurable for every $u\in U$, and
it is easily verified that $ F(\omega, s,\cdot)$
is continuous for every $(\omega,s)\in \Omega\times [0,T]$.
By a
classical selection theorem (see \cite{AuFr}, Theorems 8.1.3 and 8.2.11)
 there exists a function $u^{*,t,x}:\Omega\times [t,T]\to U$,
 measurable with respect to $\bar\calg$ and
$\calu$, such that $F(\omega,s ,u^{*,t,x}(\omega,s)) =
\min_{u\in U} F (\omega,s,u)$ for every $(\omega,s)\in \Omega\times [t,T]$, so that
 \eqref{minselector} holds true for every $(\omega,s)$.
 Note that $u^{*,t,x}$ may depend on $t,x$ because $\mu$ does.
 After modification on a set of $\mu$-measure
zero, the function $u^{*,t,x}$ can be made
measurable with respect to $\calp^t\otimes \calb([t,T])$ and
$\calu$, and \eqref{minselector} still holds, as it is understood as an equality
for $\mu$-almost
all $(\omega,s)$.
\qed

\end{document}